\documentclass[11pt]{article}

\usepackage{amsfonts,amsmath}
\usepackage{latexsym}
\usepackage{bbm}
\usepackage{todonotes}

\author{K\'aroly J. B\"or\"oczky\footnote{Supported by
NKFIH grant K 132002}\footnote{Alfr\'ed R\'enyi Institute of Mathematics, 
 Realtanoda u. 13-15, H-1053 Budapest, Hungary,
 boroczky.karoly.j@renyi.hu}}
\title{The Logarithmic Minkowski conjecture and the $L_p$-Minkowski Problem}

\newcommand{\R}{\mathbb{R}}

\newcommand{\HH}{\mathcal{H}}

\newtheorem{lemma}{LEMMA}[section]
\newtheorem{theo}[lemma]{THEOREM}

\newtheorem{conj}[lemma]{CONJECTURE}

\begin{document}

\maketitle

\begin{abstract}
The current state of art concerning the
  $L_p$-Minkowski problem as a Monge-Amp\'ere equation on the sphere and Lutwak's Logarithmic Minkowski conjecture about the uniqueness of even solution in the $p=0$ case are surveyed and connections to many related problems are discussed.
\end{abstract}

\noindent{\bf MSC 2010 } Primary: 35J96, Secondary: 52A40

\section{Introduction}

The Minkowski problem forms the core of various areas in fully nonlinear partial differential equations and convex geometry
(see Trudinger, Wang \cite{wang} and Schneider \cite{Sch14}), which was
extended to the $L_p$-Minkowski theory
 by Lutwak \cite{Lut93,Lut93a,Lut96} where $p=1$ corresponds to the classical case.
The classical Minkowski's existence theorem due to Minkowski and Aleksandrov characterizes the surface area measure
$S_K$ of a convex body $K$ in $\R^n$, more precisely, it solves the Monge-Amp\'ere equation
$$
\det(\nabla^2 h+h\,{\rm Id})= f
$$
on the sphere $S^{n-1}$ where a convex body $K$ with $C^2_+$ boundary provides a solution if
$h=h_K|_{S^{n-1}}$ for the support function $h_K$ of $K$, and in this case, $1/f(u)$ is the Gaussian curvature at the $x\in\partial K$ where $u$ is an exterior normal for $u\in S^{n-1}$. The so-called log-Minkowski problem (short hand for logarithmic Minkowski problem) 
\begin{equation}
\label{logMinkintro}
h\det(\nabla^2 h+h\,{\rm Id})= f
\end{equation}
or $L_0$-Minkowski problem was posed by Firey in his seminal paper \cite{Fir74}. It seeks to charaterize the cone volume measure $dV_K=\frac1n\,h_KdS_K$ of a convex body $K$ containing the origin $o$, and to determine whether the even solution is unique if $f$ is even. The latter problem is the so-called Log-Minkowski conjecture by Lutwak.
However, the log-Minkowski problem has received due attention only after finding its place as part of Lutwak's $L_p$-Minkowski problem 
$$
h^{1-p}\det(\nabla^2 h+h\,{\rm Id})= f
$$
 in the 1990's where the cases $p=1$ and $p=0$ are the classical and the logarithmic Minkowski problem. 

The uniqueness of the solution of the classical Minkowski problem up to translation is the consequence of the
equality case of the Brunn-Minkowski inequality (see Gardner \cite{Gar02} or Schneider \cite{Sch14})
\begin{equation}
\label{BM}
V(\alpha K+\beta C)^{\frac{1}n}\geq \alpha V(K)^{\frac{1}n}+ \beta V(C)^{\frac{1}n}
\end{equation}
for convex bodies $K$ and $C$ and $\alpha,\beta\geq 0$.  
For $p\geq 0$, the $L_p$-Minkowski problem is intimately related to the $L_p$ version of the Brunn-Minkowski inequality/conjecture
\begin{equation}
\label{lpineqconj}
\begin{array}{rcll}
V((1-\lambda)K+_p\lambda C)^{\frac{p}n}&\geq&(1-\lambda)\,V(K)^{\frac{p}n} +\lambda\,V(C)^{\frac{p}n}&\mbox{ if $p>0$;}\\[1ex]
V((1-\lambda)K+_0\lambda C)&\geq&V(K)^{1-\lambda} V(C)^\lambda & \mbox{ if $p=0$}
\end{array} 
\end{equation}
for $\lambda\in(0,1)$ and convex bodies $K,C$ containing the origin. Here \eqref{lpineqconj} is the classical Brunn-Minkowsi inequality  
if $p=1$, a theorem of Firey \cite{Fir62} if $p>1$, and assuming that $K$ and $C$ are origin symmetric, a conjecture being the central theme of this paper if $p\in[0,1)$. Actually, the conjecture has been recently verified if $p\in(0,1)$ is close to $1$; more precisely, combining Kolesnikov, Milman \cite{KoM22} and Chen, Huang, Li, Liu \cite{CHLL20} proves that
the $L_p$ Brunn-Minkowski conjecture holds if $p\in[p_n,1)$ and $K,C$ are origin symmetric convex bodies
for an explicit $p_n\in(0,1)$.

The main goal of this survey is to inspire the resolution of the Log-Brunn-Minkowski conjecture ({\it cf.} \eqref{lpineqconj} when $p=0$), or Lutwak's essentially equivalent Log-Minkowski conjecture (the Monge-Amp\'ere equation \eqref{logMinkintro} 
on $S^{n-1}$ has a unique even solution if $f$ is even, positive and $C^\infty$). 

Its versatility is an intriguing aspect 
 of the log-Minkowski conjecture; namely, uniqueness of the even solutions of a Monge-Amp\'ere equation  on the sphere is  equivalent to some strengthening of the Brunn-Minkowski inequality for origin symmetric convex bodies, to an inequality for the Gaussian density, and to some spectral  gap estimates for certain self-adjoint operators, and  in turn to displacement convexity of certain functional of probability measures on the sphere in optimal transportation.

In this survey, we review some related aspects of the classical Brunn-Minkowski Theory in Section~\ref{secMinkowski},
the state of art concerning the Log-Minkowski problem and the Log-Minkowski conjecture in
Section~\ref{secLogMinkowski}, Lutwak's $L_p$-Minkowski problem and the $L_p$-Minkowski conjecture in
Section~\ref{secLpMinkowski}, and some variants of the $L_p$-Minkowski problem in
Section~\ref{secVariants}.

\section{Classical Brunn-Minkowski Theory}
\label{secMinkowski}

This section serves as the introduction into the relevant aspects of Brunn-Minkowski Theory on the one hand, but also introduces the basic ideas and tools used in  the upcoming sections. For a thorough discussion the subject and related problems from various perpectives, see 
Artstein-Avidan, Giannopoulos, Milman \cite{AGM15,AGM21},
 Gardner \cite{Gar02}, Leichtwei{\ss} \cite{Lei98}, Schneider \cite{Sch14}.

We call a compact convex set in $\R^n$ with non-empty interior a convex body. The family of convex bodies in $\R^n$ is denoted by $\mathcal{K}^n$, and we write $\mathcal{K}_o^n$ ($\mathcal{K}_{(o)}^n$) to denote the subfamily of $K\in \mathcal{K}^n$ with $o\in K$ ($o\in{\rm int}\,K$), and write $\mathcal{K}_e^n$ to denote the family of origin symmetric convex bodies in $\R^n$. The support function of a  compact convex set $K$ is 
$h_K(u)=\max_{x\in K}\langle u,x\rangle$ for $u\in \R^n$, and hence $h_K$ is convex and homogeneous (the latter property says that
$h_K(\lambda u)=\lambda h_K(u)$ for $\lambda\geq 0$). In turn, for any convex and homogeneous function $h$ on $\R^n$, there exists a unique compact convex set $K$ such that $h=h_K$. We note that differences of support functions are dense among continous functions on the sphere; more precisely, functions of the form
$(h_K-h_C)|_{S^{n-1}}$ for convex bodies $K$ and $C$ with $C^\infty_+$ boundary in $\R^n$ are dense in $C(S^{n-1})$ with respect to the $L_\infty$ metric.

 We say that convex bodies $K$ and $C$ are homothetic  if $K=\gamma C+z$ for $\gamma>0$ and $z\in\R^n$.
We write $V(X)$ to denote Lebesgue measure of a measurable subset $X$ of $\R^n$
(with $V(\emptyset)=0$), and $\HH$ to denote the $(n-1)$-Hausdorff measure normalized in a way such that it coincides with the $(n-1)$-dimensional Lebesgue measure on ${n-1}$-dimensional affine subspaces. 
For  $X,Y\subset \R^n$ and $\alpha,\beta\in\R$, the Minkowski linear combination is
$\alpha X+\beta Y=\{\alpha x+\beta y:x\in X,\;y\in Y\}$, which is convex compact if $X$ and $Y$ are convex compact.
We write $B^n$ to denote the unit Euclidean ball centered at the origin $o$, and equip the space of compact convex sets of $\R^n$ with topology induced by the Hausdorff metric (sometimes called Hausdorff distance); namely, if $K$ and $C$ are compact convex sets, then their Haussdorff distance is
$$
\delta_H(K,C)=\min\{r\geq 0:K\subset C+rB^n\mbox{ and } C\subset K+rB^n\}.
$$

The Brunn-Minkowski inequality 
 says that if $\alpha,\beta>0$ and $K,C$ are convex bodies in $\R^n$, then
\begin{equation}
\label{BrunnMinkowski}
V(\alpha K+\beta C)^{\frac1n}\geq \alpha V(K)^{\frac1n}+\beta V(C)^{\frac1n},
\end{equation}
with equality if and only if
$K$ and $C$  are homothetic. We note that the Brunn-Minkowski inequality \eqref{BrunnMinkowski} also holds
if $K$ and $C$ are bounded Borel subsets of $\R^n$ (note that Minkowski linear combination of measurable subsets may not be 
 measurable; therefore, outer measure is used in that case).
The Brunn-Minkowski inequality famously yields the isoperimertric inequality; namely, that the surface area of a bounded Borel set $X$ of given volume is minimized by balls. Naturally, one needs a suitable notion of surface area. It is the Hausdorff measure $\HH(\partial X)$ if
$X$ is a convex body, or more generally, $\partial X$ is the finite union of the images of Lipschitz functions defined on bounded subsets of $\R^{n-1}$ (see Schneider \cite{Sch14} or Ambrosio, Colesanti, Villa \cite{ACV08}), but the right notion is finite perimeter (see Maggi \cite{Mag12}).
Fusco, Maggi, Pratelli \cite{FMP08} proved an optimal stability version of the isoperimetric inequality in terms of the symmetric difference metric (whose result was extended to the Brunn-Minkowski inequality by Figalli, Maggi, Pratelli \cite{FMP09,FMP10}, see Theorem~\ref{Maggi} below).

Because of the homogeneity of the Lebesgue measure, an equivalent form of the Brunn-Minkowski inequality \eqref{BrunnMinkowski} is that
if $K,C$ are convex bodies in $\R^n$ and $\lambda\in(0,1)$, then
\begin{equation}
\label{BrunnMinkowskiprod}
V((1-\lambda)K+\lambda C)\geq V(K)^{1-\lambda}V(C)^{\lambda},
\end{equation}
with equality if and only if $K$ and $C$ are translates. A big advantage of this product form of the Brunn-Minkowski inequality
is that it is dimension invariant.

The first stability forms of the Brunn-Minkowski inequality were due to Minkowski himself 
(see Groemer \cite{Gro93}).
If the distance of the convex bodies $K$ and $C$ is measured in terms of the so-called Hausdorff distance, then
Diskant \cite{Dis73} and Groemer \cite{Gro88} provided close to optimal stability versions
(see Groemer \cite{Gro93}). However, the natural distance is in terms of the volume of the symmetric difference, and the optimal result is due to Figalli, Maggi, Pratelli \cite{FMP09,FMP10}.
To define the ``homothetic distance'' $A(K,C)$
of convex bodies $K$ and $C$, let $\alpha=|K|^{\frac{-1}n}$ and
$\beta=V(C)^{\frac{-1}n}$, and let
$$
A(K,C)=\min\left\{V(\alpha K\Delta (x+\beta C)):\,x\in\R^n\right\}.
$$
In addition, let
$$
\sigma(K,C)=\max\left\{\frac{V(C)}{V(K)},\frac{V(K)}{V(C)}\right\}.
$$

\begin{theo}[Figalli, Maggi, Pratelli]
\label{Maggi}
For $\gamma^*(n)>0$ depending on $n$
and  any convex bodies $K$ and $C$ in $\R^n$,
$$
V(K+C)^{\frac1n}\geq (V(K)^{\frac1n}+V(C)^{\frac1n})
\left[1+\frac{\gamma^*(n)}{\sigma(K,C)^{\frac1n}}\cdot A(K,C)^2\right].
$$
\end{theo}

Here the exponent $2$ of $A(K,C)^2$ is optimal, see Figalli, Maggi, Pratelli \cite{FMP10}.
We note that prior to \cite{FMP10}, the only known error term in the Brunn-Minkowski inequality
was of order $A(K,C)^\eta$ with $\eta\geq n$,
coming from the estimates of 
Diskant \cite{Dis73} and  Groemer \cite{Gro88} in terms of the Hausdorff distance.

Figalli, Maggi, Pratelli \cite{FMP10} proved a factor of the form $\gamma^*(n)=cn^{-14}$ for some absolute constant $c>0$, which was improved to $c n^{-7}$ by Segal \cite{Seg12}, and subsequently to $c n^{-5.5}$ by Kolesnikov, Milman \cite{KoM22}, Theorem~12.12. 
The current best known bound for $\gamma^*(n)$ is $cn^{-5}(\log n)^{-10}$, which follows by combining the general estimate of Kolesnikov-Milman \cite{KoM22}, Theorem 12.2, with the polylogarithmic bound of Klartag, Lehec \cite{KlL} on the Cheeger constant of a convex body in isotropic position improving on Yuansi Chen’s work \cite{Che21} on the Kannan-Lov\'asz-Simonovits conjecture.
Harutyunyan \cite{Har18} conjectured that $\gamma^*(n)=c n^{-2}$ is the optimal order of the constant, and 
showed that it can't be of smaller order. Actually, Segal \cite{Seg12} observed that
 Dar's conjecture in \cite{Dar99} would imply that we may choose $\gamma^*(n)=c n^{-2}$
for some absolute constant $c>0$. Here Dar's conjectured strengthening of the Brunn-Minkowski inequality states in \cite{Dar99} that if $K$ and $C$ are convex bodies in $\R^n$, and
$M=\max_{x\in\R^n}V(K\cap (x+C))$, then
\begin{equation}
\label{Dar}
V(K+C)^{\frac1n}\geq M^{\frac1n}+\left(\frac{V(K)V(C)}{M}\right)^{\frac1n}.
\end{equation}
Dar's conjecture is only known to hold in the plane (see Xi, Leng \cite{XiL16}), and in some very specific cases in higher dimension (see Dar \cite{Dar99}). 

The paper Eldan, Klartag \cite{ElK14} discusses "isomorphic" stability versions of the Brunn-Minkowski inequality under condition of the type
$|\frac12\,K+\frac12\,C|\leq 5\sqrt{|K|\cdot|C|}$,
and considers, for example, the $L^2$ Wasserstein distance of the uniform measures on suitable affine images of $K$ and $C$.

We note that stability versions of the Brunn-Minkowski inequality have been verified even if $K$ or $C$ are not convex.
The essentially optimal estimate Theorem~\ref{Maggi} (with much worse factor $\gamma^*(n)$) 
is verified if $K$ is bounded measurable and $C$ a convex body by Barchiesi, Julin \cite{BaJ17}
(improving on the estimate in Carlen, Maggi \cite{CaM17}), if $n\geq 2$ and $K=C$ is bounded Borel set by
Hintum, Spink, Tiba \cite{HST}, and if $n=2$ and $K$ and $C$ are bounded Borel sets by
Hintum, Spink, Tiba \cite{HST2}. If $n\geq 3$ and $K$ and $C$ are bounded Borel sets, then only a much weaker estimate
in terms of $A(K,C)$ is known, proved by Figalli, Jerison \cite{FiJ15,FiJ17}.
On the other hand, a better error term  of order $A(X,Y)$ holds if $n=1$  according to Freiman
and Christ (see Christ \cite{Chr12}). 

It was proved by Minkowski that if $K$ and $C$ are convex bodies and $\alpha,\beta\geq 0$, then
\begin{equation}
\label{mixedvolume}
V(\alpha\,K+\beta\,C)=\sum_{i=0}^n{n \choose i}V(K,C;i)\alpha^{n-i}\beta^i
\end{equation}
where $V(K,C;i)$ are the so-called mixed volumes. For fixed $i$, $V(K,C;i)$ is positive, continuous in both variables, and satisfies $V(\alpha K,\beta C;i)=\alpha^{n-i}\beta^iV(K,C;i)$ for $\alpha,\beta>0$, 
$V(K,C;i)=V(C,K;n-i)$ and $V(\Phi K+x,\Phi C+y;i)=V(K,C;i)$ for $x,y\in\R^n$ and $\Phi\in{\rm SL}(n)$.
Many mixed volumes have geometric meaning, like $V(K,K;i)=V(K,C;0)=V(K)$ and 
$\frac1n\,V(K,B^n;1)=\HH(\partial K)$ is 
the surface area of $K$. In addition, if $i=1,\ldots,n-1$, then $V(K,B^n;n-i)$ is proportional with the mean $i$-dimensional projection of $K$ according to the Kubota formula (see Leichtwei{\ss} \cite{Lei98}, Schneider \cite{Sch14}).

It follows from the Brunn-Minkowski inequality \eqref{BrunnMinkowski} that the function
$f(\lambda)=V((1-\lambda)K+\lambda C)^{\frac1n}$ is concave on $[0,1]$. Combining $f'(0)\geq f(1)-f(0)$ and 
\eqref{mixedvolume} leads to the famous Minkowski inequality
\begin{equation}
\label{MinkowskiIneq}
V(K,C;1)^n\geq V(K)^{n-1}V(C),
\end{equation}
with equality if and only if $K$ and $C$ are homothetic. The Minkowski inequality is actually equivalent to the Brunn-Minkowski inequality because it implies that the function $f(\lambda)=V((1-\lambda)K+\lambda C)^{\frac1n}$ is concave on $[0,1]$, which in turn yields \eqref{BrunnMinkowski}. Considering the second derivative $f"(\lambda)$ leads to Minkowski's second inequality
\begin{equation}
\label{Minkowski2ineq}
V(K,C;1)^2\geq V(K)V(K,C;2),
\end{equation}
that is in turn also equivalent to the Brunn-Minkowski inequality \eqref{BrunnMinkowski}. The rather involved equality case of \eqref{Minkowski2ineq} has been only recently clarified by van Handel, Shenfeld \cite{vHS22}.

Actually, Minkowski defined the mixed volume $V(C_1,\ldots,C_n)$ of $n$ convex bodies  via the identity
\begin{equation}
\label{Mixedvol}
V\left(\lambda_1K_1+\ldots+\lambda_mK_m\right)=\sum_{i_1,\ldots,i_n=1}^m V(K_{i_1},\ldots,K_{i_n})
\cdot \lambda_{i_1}\cdot\ldots\cdot\lambda_{i_n}
\end{equation}
for $K_1,\ldots,K_m\in\mathcal{K}^n$ and $\lambda_1,\ldots,\lambda_m\geq 0$ where 
$V(C_1,\ldots,C_n)\geq 0$ is symmetric and continuous in its variables  (see \cite{Sch14}), and $V(K,C;i)$ means $i$ copies of $C$ and $n-i$ copies of $K$. A far reaching generalization of Minkowski's first and second  inequalities is the Alexandrov-Fenchel inequality
$$
V(K_1,K_2,K_3,\ldots,K_n)^2\geq V(K_1,K_1,K_3,\ldots,K_n)V(K_2,K_2,K_3,\ldots,K_n)
$$
(see  Alexandrov \cite{Ale37a,Ale96} and Schneider \cite{Sch14}). Equality in the Alexander-Fenchel inequality has been much better understood now due to van Handel, Shenfeld \cite{vHS22,vHS} where \cite{vHS} clarifies the case of polytopes. 

 Let us summarize some equivalent formulations of the Brunn-Minkowski inequality that holds for all convex bodies $K$ and $C$ in 
$\R^n$:
\begin{itemize}
\item $V(\alpha K+\beta C)^{\frac1n}\geq \alpha V(K)^{\frac1n}+\beta V(C)^{\frac1n}$ (({\it cf.} \eqref{BrunnMinkowski});

\item $V((1-\lambda)K+\lambda C)\geq V(K)^{1-\lambda}V(C)^{\lambda}$ (({\it cf.} \eqref{BrunnMinkowskiprod});

\item $f(\lambda)=V((1-\lambda)K+\lambda C)^{\frac1n}$ is concave on $[0,1]$;

\item Minkowski inequality
$V(K,C;1)^n\geq V(K)^{n-1}V(C)$ ({\it cf.}  \eqref{MinkowskiIneq});

\item Minkowski's second inequality $V(K,C;1)^2\geq V(K)V(K,C;2)$ ({\it cf.} \eqref{Minkowski2ineq}).

\end{itemize}

The classical Minkowski problem  is concerned with the characterization of the so-called surface area measure $S_K$ of a convex body $K$. Let $\partial'K$ denote the subset of the boundary  of $K$ such that there exists a unique exterior unit normal vector
$\nu_K(x)$ at any point $x\in \partial'K$.
It is well-known that $\HH(\partial K\setminus\partial'K)=0$ and $\partial'K$ is a Borel set  (see Schneider \cite{Sch14}). The function $\nu_K:\partial'K\to S^{n-1}$ is the spherical Gauss map
that is continuous on $\partial'K$.
 The surface area measure $S_K$ of $K$ is a Borel measure on $S^{n-1}$ 
satisfying that $S_K(\eta)=\HH(\nu_K^{-1}(\eta))$
 for any Borel set $\eta\subset S^{n-1}$. The surface area measure is the first variation of the volume; namely,
if $C$ is any convex body in $\R^n$, then
\begin{equation}
\label{SurfaceAreaVolume}
nV(K,C;1)=\lim_{\varepsilon\to 0^+}\frac{V(K+\varepsilon C)-V(K)}{\varepsilon}=\int_{S^{n-1}}h_C\,d S_K.
\end{equation}
It also follows that the Minkowski inequality \eqref{MinkowskiIneq} can be written in the following form: If 
$V(K)=V(C)$ holds for $K,C\in\mathcal{K}^n$, then
\begin{equation}
\label{MinkIneqSK}
\int_{S^{n-1}}h_C\,dS_K\geq 
\int_{S^{n-1}}h_K\,dS_K,
\end{equation}
with equality if and only if $K$ and $C$ are translates.  

To consider some examples for the surface area measure, if $P$ is a polytope with facets $F_1,\ldots,F_k$ and exterior unit normals $u_1,\ldots, u_k$, then $S_P$ is concentrated onto 
$\{u_1,\ldots, u_k\}$ and $S_P(u_i)=\HH(F_i)$ for $i=1,\ldots,k$. On the other hand, if $\partial K$ is $C^2_+$, namely, $C^2$ with positive Gaussian curvature, and the Gaussian curvature at the point of $\partial K$ with exterior unit normal $u\in S^{n-1}$ is $\kappa(u)=\kappa(K,u)$, then 
\begin{equation}
\label{SKGauss}
dS_K=\kappa^{-1}\,d\HH=\det(\nabla^2 h+h\,{\rm Id})\,d\HH
\end{equation}
on $S^{n-1}$ where $h=h_K|_{S^{n-1}}$ and $\nabla h$ and  $\nabla^2 h$ are the gradient and the Hessian of $h$ with respect to a moving orthonormal frame. In particular, $S_K$ is absolute continuous in this case.

We note that if $\partial K$ is $C^2_+$ for $K\in\mathcal{K}^n$ and $h=h_K|_{S^{n-1}}$, then for any $u\in S^{n-1}$, the differential operator
\begin{equation}
\label{D2h}
D^2h(u)=\nabla^2 h(u)+h(u)\,{\rm Id}
\end{equation}
is the restriction of the Hessian of $h_K$ (in $\R^n$) at $\lambda u$ to an operator $u^\bot\mapsto u^\bot$ for any $\lambda>0$, and the eigenvalues of $D^2h(u)$ are the radii of curvature at $x\in \partial K$ where $u$ the exterior unit normal is. In turn, for any given $h\in C^m(S^{n-1})$ with $m\geq 2$, $h=h_K|_{S^{n-1}}$  for $K\in\mathcal{K}^n$ with $C^m$ ($C^m_+$) boundary if and only if $D^2h(u)$ is positive semi-definit (positive definit) for $u\in S^{n-1}$.

Now the Minkowski problem asks for necessary and sufficient conditions for a Borel measure $\mu$ on $S^{n-1}$
such that
\begin{equation}
\label{Minkprobmu}
\mu=S_K
\end{equation}
for a convex body $K$. The solution, together with its uniqueness was provided by Minkowski \cite{Min97,Min03}
if the measure $\mu$ is discrete (and hence the convex body is a polytope) or absolutely continuous.
Minkowski's solution was extended to any general measure $\mu$ by Alexandrov \cite{Ale37b,Ale38,Ale96}; namely,
there exists a convex body $K$ with $\mu=S_K$ if and only if
\begin{eqnarray}
\label{Minkprob1}
\mu(L\cap S^{n-1})<\mu(S^{n-1})&&\mbox{for any linear $(n-1)$-subspace $L\subset \R^n$;}\\
\label{Minkprob2}
\int_{S^{n-1}}u\,dS_K(u)=o;&&
\end{eqnarray}
moreover, $S_K=S_C$ holds for convex bodies $K$ and $C$  if and only if $K$ and $C$ are translates.
Essentially complete solutions were published also by Fenchel, Jensen \cite{FeJ38} and Lewy \cite{LE} about the same time. 
In particular, the Monge-Amp\`ere equation on the sphere $S^{n-1}$ corresponding to the Minkowski problem is
\begin{equation}
\label{MongeSK}
\det(\nabla^2 h+h\,{\rm Id})= f
\end{equation}
where $f$ is a given non-negative function with positive integral.
If the given Borel measure $\mu$ on $S^{n-1}$ is not absolutely continuous with respect to the Lebesgue measure, then
$h=h_K|_{S^{n-1}}$ is a solution of \eqref{MongeSK} in the Alexandrov sense if \eqref{Minkprobmu} holds.
We note that the surface area measure $S_K$ is actually the  Monge-Amp\`ere measure corresponding to $h$ (see Trudinger, Wang \cite{wang} and B\"or\"oczky, Fodor \cite{BoF19}, Section~7).

The regularity of the solution of the Minkowski problem \eqref{MongeSK} is well investigated by
 Nirenberg \cite{NIR}, Cheng and Yau \cite{CY}, Pogorelov \cite{POG}, showing eventually that if $f$ is positive and $C^k$ for $k\geq 1$, then $h$ is $C^{k+2}$. Finally, Caffarelli \cite{Caf90a,Caf90b} proves that if $f$ is positive and $C^\alpha$ for $\alpha\in(0,1)$ (namely, $|f(x)-f(y)|\leq C\|x-y\|^{\alpha}$ for $x,y\in S^{n-1}$ and constant $C>0$), then the solution $h$ is $C^{2,\alpha}$ (see also B\"or\"oczky, Fodor \cite{BoF19}, Section~7 on how to connect results on Monge-Amp\`ere equations on $\R^n$ to Monge-Amp\`ere equations on $S^{n-1}$, and
Chen, Liu, Wang \cite{CLW21} for an extension of \cite{Caf90a,Caf90b}).

Turning to proofs, one of the elegant arguments proving the Brunn-Minkowski inequality \eqref{BrunnMinkowski} is due to Hilbert, and is based on a spectral gap estimate for a differential operator (see \eqref{Hilbert-operator} and Bonnesen, Fenchel \cite{BoF87}). This approach was further developed by
Alexandrov \cite{Ale37a,Ale96} leading to the Alexandrov-Fenchel inequality, by van Handel, Shenfeld \cite{vHS22,vHS} to characterize equality in the Alexandrov-Fenchel inequality in certain case, and by
Milman, Kolesnikov \cite{KoM22} leading to the $L_p$-Minkowski inequlity Theorem~\ref{pcloseto1} improving the
Brunn-Minkowski inequality for origin symmetric convex bodies (see the end of Section~\ref{secLpMinkowski}).
Another fundamental approach proving the Brunn-Minkowski inequality is initiated by
Gromov's influential appendix to Milman, Schechtman \cite{MSG86} using ideas by
 Knothe \cite{Kno57} provided a proof of the isoperimetric inequality using optimal (mass) transport, and the argument can be readily extended to the Brunn-Minkowski inequality  \eqref{BrunnMinkowski} and the Pr\'ekopa-Leindler inequality  \eqref{PLnineq} below.
This approach lead even to the stability version Theorem~\ref{Maggi} by Figalli, Maggi, Pratelli \cite{FMP09,FMP10}.
We note that the original argument of Brunn and Minkowski for \eqref{BrunnMinkowski} (see Bonnesen, Fenchel \cite{BoF87})  can be also considered as a version of the mass transportation approach.

For the Minkowski problem \eqref{Minkprobmu}, the variational method seeks the minimum of
$\int_{S^{n-1}}h_C\,d\mu$ over all convex bodies $C$ with $V(C)=1$ where $\mu$ satisfies \eqref{Minkprob1} and
\eqref{Minkprob2}. It follows from \eqref{Minkprob2} that the integral is invariant under translating $C$, hence the existence of a minimizer $C_0$ can be established. The fact that $S_{C_0}$ is proprotional to $\mu$ follows via Alexandrov's Lemma~\ref{AlexandrovLemma} extending \eqref{SurfaceAreaVolume} (see Theorem~7.5.3 in Schneider \cite{Sch14}).

\begin{lemma}[Alexandrov]
\label{AlexandrovLemma}
Given a convex body $K$ in $\R^n$ and continous functions $h_t,g:S^{n-1}\to\R$, 
let us assume that the Wulff shape $K_t=\{x\in\R^n:\langle x,u\rangle\leq h_t(u)\;\forall u\in S^{n-1}\}$ is a convex body and
$\lim_{t\to 0}\frac{h_t(u)-h_K(u)}t=g(u)$ uniformly in $u\in S^{n-1}$. Then
$$
\lim_{t\to 0}\frac{V(K_t)-V(K)}t=\int_{S^{n-1}}g\,dS_K.
$$
\end{lemma}

 The classical functional form of the Brunn-Minkowski inequality is the Pr\'ekopa-Leindler due to
Pr\'ekopa \cite{Pre71} and Leindler \cite{Lei72} in dimension one, was generalized in
 Pr\'ekopa \cite{Pre73,Pre75} and  Borell \cite{Bor75} ({\it cf.} also Marsiglietti \cite{Mar17}, Bueno, Pivovarov \cite{BuP}), Brascamp, Lieb \cite{BrL76}, Kolesnikov, Werner \cite{KoW22}, 
Bobkov, Colesanti, Fragal\`a \cite{BCF14}).
Various applications are provided and surveyed in Ball \cite{Bal},
Barthe \cite{Bar}, 
Fradelizi, Meyer \cite{FrM07} and Gardner \cite{Gar02}.
The following multiplicative version from \cite{Bal} is often more useful
and is more convenient for geometric applications.

\begin{theo}[Pr\'ekopa-Leindler]
\label{PLn}
If $\lambda\in(0,1)$ and $h,f,g$ are non-negative integrable functions on $\R^n$
satisfying $h((1-\lambda)x+\lambda y)\geq f(x)^{1-\lambda}g(y)^\lambda$ for
$x,y\in\mathbb{R}^n$, then
\begin{equation}
\label{PLnineq}
\int_{\R^n}  h\geq \left(\int_{\R^n}f\right)^{1-\lambda} \cdot \left(\int_{\R^n}g\right)^\lambda.
\end{equation}
\end{theo}

For a convex function $W:\R^n\to (-\infty,\infty]$, we say that the function $\varphi=e^{-W}$ is log-concave where $e^{-\infty}=0$ (in other words, $\log \varphi$ is concave for $\varphi:\R^n\to [0,\infty)$). According to Dubuc \cite{Dub77}, if equality holds in \eqref{PLn} assuming $\int_{\R^n}h>0$, then $h$ is log-concave, and
there exist $a>0$ and $z\in \R^n$ such that 
$f(x)= a^\lambda\,h(x-\lambda z)$ and $g(x)=a^{-(1-\lambda)}h(x+(1-\lambda) z)$ for almost all $x\in\R^n$. Stability versions of the Pr\'ekopa-Leindler inequality in terms of the $L_1$ distance have been established by
B\"or\"oczky, De \cite{BoD21a}  in the log-concave case, and by
 B\"or\"oczky, Figalli, Ramos \cite{BFR} for any functions where the case of log-concave functions in one variable
have been dealt with ealier by  Ball, B\"or\"oczky \cite{BaB10}. A stability version of the Pr\'ekopa-Leindler inequality
of somewhat different nature is due to  
Bucur, Fragal\`a \cite{BuF14}.

An "isomorphic" stability result for the Pr\'ekopa-Leindler inequality, in terms of the transportation distance is obtained in Eldan \cite{Eld13}, Lemma~5.2. By rather standard considerations, one can show that non-isomorphic stability results in terms of transportation distance imply stability in terms of $L_1$ distance (e.g., such implication is attained by combining Proposition~2.9 in Bubeck, Eldan, Lehec \cite{BEL18} and Proposition~10 in Eldan, Klartag \cite{ElK14}). However, the current result in \cite{Eld13}, due to its isomorphic nature, falls short of being able to obtain a meaningful bound in terms of the $L_1$ distance.

Brascamp, Lieb \cite{BrL76} proved a local version of the Pr\'ekopa-Leindler inequality for log-concave functions (Theorem~4.2 in \cite{BrL76}), which is equivalent to a Poincare-type so called Brascamp-Lieb inequality Theorem~4.1 in \cite{BrL76}. 
The paper Livshyts \cite{Liva} provides a stability version of this Brascamp-Lieb inequality, and
Bolley, Cordero-Erausquin, Fujita, Gentil, Guillin \cite{BCFGG20} proves a more general inequality.

Our final topic in this section is the Blasshke-Satal\'o inequality \eqref{Blaschke-Santalo}. 
For $K\in\mathcal{K}_{(o)}^n$ and $u\in S^{n-1}$, the radial function $\varrho_K(u)>0$ satisfies $\varrho_K(u)u\in\partial K$, and the polar (dual) $K^*\in\mathcal{K}_{(o)}^n$ of $K$ is defined by 
$\varrho_{K^*}(u)=h_K(u)$ for $u\in S^{n-1}$. 
Next, the centroid of a convex body $K$ in $\R^n$
is $\sigma_K=\frac1{V(K)}\int_Kx\,dx$, which is invariant under affine transformations. For $K\in\mathcal{K}^n$, we call it centered if $\sigma_K=o$, and Kannan, Lov\'asz, Simonovits \cite{KLS95} prove that there exists a centered ellipsoid $E$
such that $E\subset K\subset nE$ in this case.

 According to the Blashke-Santal\'o inequality (see Santal\'o \cite{San49}, Luwak \cite{Lut93} or
Schneider \cite{Sch14}), if  $K\in\mathcal{K}^n$ is centered, then $V(K^*)V(K)\leq V(B^n)^2$, or equivalently,
\begin{equation}
\label{Blaschke-Santalo}
\int_{S^{n-1}}h_K^{-n}\,d\HH\leq \frac{nV(B^n)^2}{V(K)},
\end{equation}
with equality if and only if $K$ is a centered ellipsoid.

The Blaschke-Santal\'o inequality can be proved for example via the Brunn-Minkowski inequality
(see Ball \cite{Bal} in the origin symmetric case, and Meyer, Pajor \cite{MeP90} in general). Various equivalent formulations are discussed in the beautiful survey Lutwak \cite{Lut93} (see also B\"or\"oczky \cite{Bor10} and Schneider \cite{Sch14}).
Stability versions of the Blaschke-Santal\'o inequality  are verified in 
B\"or\"oczky \cite{Bor10} and Ball, B\"or\"oczky \cite{BaB11}.
Following Ball \cite{Bal}, functional versions of the Blaschke-Santal\'o inequality have been obtained by
Artstein-Avidan, Klartag, Milman \cite{AKM04}, 
Fradelizi, Meyer \cite{FrM07}, Lehec \cite{Leh09a,Leh09b}, 
Kolesnikov, Werner \cite{KoW22},  Kalantzopoulos, Saroglou \cite{KaS}.

Recently, various breakthrough stability results about geometric functional inequalities have been obtained. 
 Stonger versions of the 
functional Blaschke-Santal\'o inequality is provided by 
Barthe, B\"or\"oczky, Fradelizi \cite{BBF14}, of the Borell-Brascamp-Lieb inequality is provided by Ghilli, Salani \cite{GhS17}, Rossi, Salani \cite{RoS17,RoS19} and Balogh, Krist\'aly \cite{BaK18} (later even on Riemannian manifolds),
 of the Sobolev inequality by Figalli, Zhang \cite{FiZ} (extending Bianchi, Egnell \cite{BiE91} and Figalli, Neumayer \cite{FiN19}), 
Nguyen \cite{Ngu16} and Wang \cite{Wan16}, of the log-Sobolev inequality by Gozlan \cite{Goz},
and of some related inequalities by 
Caglar, Werner \cite{CaW17}, Cordero-Erausquin \cite{Cor17}, Kolesnikov, Kosov \cite{KoK17}.
Another functional version of the Brunn-Minkowski inequality is provided by 
Artstein-Avidan,  Florentin, Segal \cite{AFS20}.

\section{Cone volume measure, log-Minkowski problem, log-Brunn-Minkowski conjecture}
\label{secLogMinkowski}

Given a convex body $K$ containing the origin, the cone volume measure is defined as $dV_K=\frac1n\,h_K\,dS_K$, and hence 
the total measure is $V_K(S^{n-1})=V(K)$. The name originates from the fact that if $P$ is a polytope with facets 
$F_1,\ldots,F_k$ and exterior unit normals $u_1,\ldots, u_k$, then $V_P$ is concentrated onto 
$\{u_1,\ldots, u_k\}$ and $V_P(u_i)=\frac{h_P(u_i)}n\cdot\HH(F_i)$ is the volume of the cone ${\rm conv}\{o,F_i\}$ for $i=1,\ldots,k$. 
We note that the Monge-Amp\`ere equation on the sphere $S^{n-1}$ corresponding to the logarithmic Minkowski problem is
\begin{equation}
\label{MongeVK}
h\det(\nabla^2 h+h\,{\rm Id})=  nf
\end{equation}
for a non-negative measurable function $f$ on $S^{n-1}$ with $0<\int_{S^{n-1}}f\,d\HH<\infty$. 
It follows {\it via} Caffarelli \cite{Caf90a,Caf90b} that if $f$ is positive and $C^\alpha$ for $\alpha\in(0,1)$, then the solution $h$ is $C^{2,\alpha}$, and if $f$ is positive and $C^k$ for integer $k\geq 1$, then the solution $h$ is $C^{k+2}$.
For a finite non-trivial Borel measure $\mu$ on $S^{n-1}$, a non negative function $h$ on $S^{n-1}$ that is the restriction of the support function $h_K$ for a convex body $K$ is the solution of \eqref{MongeVK} in the Alexandrov sense if 
\begin{equation}
\label{AlexandrovVK}
d\mu=dV_K=\mbox{$\frac{1}n$}\,h_K\,dS_K.
\end{equation}

A characteristic feature of the cone volume measure is that it intertwines with linear transformations; more precisely,
$V_{\Phi K}=|\det \Phi|\cdot (\Phi^{-t})_*V_K$
for any $K\in\mathcal{K}_o^n$ and $\Phi\in{\rm GL}(n,\R)$. We note that if $u\in S^{n-1}$ is an exterior normal at an $x\in \partial K$, then $\Phi^{-t}u$ is an exterior normal at $\Phi x\in\partial(\Phi K)$, and 
 the push forward measure $\Psi_*\mu$ on $S^{n-1}$ for 
  a Borel measure $\mu$ on $S^{n-1}$ and $\Psi\in{\rm GL}(n)$  is defined (with a slight abuse of notation) in a way such that
if $\omega\subset S^{n-1}$ is Borel, then
$$
\Psi_*\mu(\omega)=\mu\left(\left\{\frac{\Psi^{-1}(u)}{\|\Psi^{-1}(u)\|}:\,u\in\omega\right\}\right).
$$

Cone volume measure was introduced by Firey \cite{Fir74}, and has been a widely used tool since the paper
Gromov, Milman \cite{GromovMilman}, see for example Barthe, Gu\'{e}don, Mendelson, Naor \cite{BG}, Naor \cite{Nar07},  Paouris, Werner \cite{PaW12}.
The still open logarithmic Minkowski problem \eqref{AlexandrovVK} or \eqref{MongeVK} was posed by Firey \cite{Fir74} in 1974, who showed that
if $f$ is a positive constant function, then \eqref{MongeVK} has a unique even solution coming from the suitable centered ball. 
For a positive constant function $f$, the general uniqueness result without the evenness condition is due to  Andrews \cite{And99} if $n=2,3$, and Brendle, Choi, Daskalopoulos \cite{BCD17} if $n\geq 4$. 
It is known that uniqueness of the solution may not hold if $f$ is not a constant function (see, for example, Chen, Li, Zhu \cite{CLZ19}).
However, the celebrated "Logarithmic Minkowski conjecture" by 
 Lutwak \cite{Lut93a} from 1993 states that
\eqref{MongeVK} has a unique even solution if $f$ is even and positive (Conjecture~\ref{logMinkowskiUniqCinfty} is a more restricted version).  

\begin{conj}[Log-Minkowski conjecture \#1]
\label{logMinkowskiUniqCinfty}
If $f$ is a positive even $C^\infty$ function in \eqref{MongeVK}, then \eqref{MongeVK} has a unique even solution.
\end{conj}

As we explain below, the following logarithmic analogue of Minkowski's inequality \eqref{MinkIneqSK} is an intimately related form of the  Log-Minkowski conjecture (see B\"or\"oczky, Lutwak, Yang, Zhang \cite{BLYZ12} for origin symmetric bodies, and 
by B\"or\"oczky, Kalantzopoulos \cite{BoK22} for centered convex bodies).

\begin{conj}[Log-Minkowski conjecture \#2]
\label{logMconj}
If $K$ and $C$ are convex bodies in $\R^n$ whose centroid is the origin, then
\begin{equation}
\label{logMconjeq}
\int_{S^{n-1}}\log \frac{h_C}{h_K}\,dV_K\geq \frac{V(K)}n\log\frac{V(C)}{V(K)}
\end{equation}
with equality if and only if $K=K_1+\ldots + K_m$ and $C=C_1+\ldots + C_m$ for compact convex sets
	$K_1,\ldots, K_m,C_1,\ldots,C_m$ of dimension at least one where $\sum_{i=1}^m{\rm dim}\,K_i=n$
	and $K_i$ and $C_i$ are dilates, $i=1,\ldots,m$.
\end{conj}

In particular, the more precise form of the Logarithmic Minkowski Conjecture~\ref{logMinkowskiUniqCinfty} is that
if $K,C\in\mathcal{K}^n$ are centered, then $V_K=V_C$ implies the equality conditions in  Conjecture~\ref{logMconj}.

We note that the choice of the right translates of $K$ and $C$ are important in Conjecture~\ref{logMconj} according to the examples by Nayar, Tkocz \cite{NaT13},
and that Conjecture~\ref{logMconj} is invariant under applying the same non-singular linear transformation to $C$ and $K$. An equivalent form of Conjecture~\ref{logMconj} is that if $V(C)=V(K)$ for centered $C$ and $K$, then  
\begin{equation}
\label{logMconjeqVCK}
\int_{S^{n-1}}\log h_C\,dV_K\geq \int_{S^{n-1}}\log h_K\,dV_K
\end{equation}
where the case of equality is like in Conjecture~\ref{logMconj}.

Let me explain why Conjecture~\ref{logMinkowskiUniqCinfty} is equivalent to the case of Conjecture~\ref{logMconj} when
$K\in\mathcal{K}_e^n$ has $C^\infty_+$ boundary. Since $V_K$ satisfies the strict subspace concentration condition (see below), B\"or\"oczky, Lutwak, Yang, Zhang  \cite{BLYZ13} prove that the function $C\mapsto \int_{S^{n-1}}\log h_C\,dV_K$ of origin symmetric convex bodies $C$ with $V(C)=V(K)$ attains its minimum; moreover, whenever it attains its minimum at some $C=\widetilde{C}$, then $V_{\widetilde{C}}=V_K$. In turn, the stated equivalence follows. In addition, it follows by approximation that Conjecture~\ref{logMinkowskiUniqCinfty} yields the inequality
\eqref{logMconjeq} for any pair of origin symmetric convex bodies $K$ and $C$ without the case of equality.

In $\R^2$, Conjecture~\ref{logMconj} is verified in B\"or\"oczky, Lutwak, Yang, Zhang \cite{BLYZ12} for origin symmetric convex bodies, but it is still open in general even in the plane.
In higher dimensions, Conjecture~\ref{logMconj} is proved for complex bodies ({\it cf.} Rotem \cite{Rotem}), and if there exist $n$ independent linear reflections that are common symmetries of $K$ and $C$ ({\it cf.} B\"or\"oczky, Kalantzopoulos \cite{BoK22}, and even a stability version is verified by B\"or\"oczky and De \cite{BoD}). The latter type of bodies include unconditional convex bodies, which case was handled earlier by Saroglou \cite{Sar15}.
In addition, Conjecture~\ref{logMconj} is verified if $C$ is origin symmetric and $K$ is a zonoid by van Handel \cite{vHa} (with equality case only clarified when $K$ has $C^2_+$ boundary), or if $C$ is a centered convex body and $K$ is a centered ellipsoid by Guan, Ni \cite{GuN17}. The latter case directly follows from the Jensen inequality and the Blaschke-Santal\'o inequality \eqref{Blaschke-Santalo}, as assuming that $K=B^n$ and $V(C)=V(B^n)$, we have
\begin{align*}
\exp\left(\int_{S^{n-1}}\log h_C\cdot \frac1{V(B^n)}dV_K\right)&=
\exp\left(\int_{S^{n-1}}\log h_C\cdot \frac1{nV(B^n)}\,d\HH\right)\\
&\geq \left(\int_{S^{n-1}}h_C^{-n}\cdot \frac1{nV(B^n)}\,d\HH\right)^{\frac{-1}n}
\geq 1.
\end{align*}

For origin symmetric $K$ and $C$, Conjecture~\ref{logMconj} is proved when $K$ is close to be an ellipsoid (with equality case only clarified when $K$ has $C^2_+$ boundary)  by a combination of the local estimates by Kolesnikov, Milman \cite{KoM22}, and the use of the continuity method in PDE by Chen, Huang, Li, Liu \cite{CHLL20}. Here closeness to an ellipsoid means that there exist some $c_n>0$ depending only on $n$ and an origin symmetric ellipsoid $E$ such that $E\subset K\subset (1+c_n)E$. Another even more recent proof of this result is due to Putterman \cite{Put21}. We note that an analogues result holds for linear images of Hausdorff neighbourhoods of $l_q$ balls for $q>2$ if the dimension $n$ is high enough according to \cite{KoM22} and the method of \cite{CHLL20}.
Actually, Milman \cite{Milb} provides rather generous explicit curvature pinching bounds for $\partial K$ in order to Conjecture~\ref{logMconj} to hold, and proves that for any origin symmetric convex body $M$ there exists an
origin symmetric convex body $K$ with $C^\infty_+$ boundary and $M\subset K\subset 8M$ such that Conjecture~\ref{logMconj} holds for any origin symmetric convex body $C$.
Additional local versions of Conjecture~\ref{logMconj} are due to Colesanti, Livshyts, Marsiglietti \cite{CLM17}, Kolesnikov, Livshyts \cite{KoL22} and Hosle, Kolesnikov, Livshyts \cite{HKL21}. We review Kolesnikov and Milman's approach in \cite{KoM22} based on the Hilbert-Brunn-Minkowski operator at the end of Section~\ref{secLpMinkowski}.

Xi, Leng \cite{XiL16} considered a version of Conjecture~\ref{logMconj} where the convex bodies $K$ and $C$ in $\R^n$ are translated by vectors depending in both $K$ and $C$. We set $r(K,C)=\max\{t>0:\exists x,\;x+tC\subset K\}$ and
$R(K,C)=\min\{t>0:\exists x,\;K\subset x+tC\}$, and say that $K$ and $C$ are in dilated position if $o\in K\cap C$ and
\begin{equation}
\label{dilated}
r(K,C)\,C\subset K\subset R(K,C)\,C.
\end{equation}
We observe that $r(C,K)\,K\subset C\subset R(C,K)\,K$ in this case. Now for any convex bodies $K$ and $C$ there exist $z\in K$ and $w\in C$ such that $K-z$ and $C-w$ are in dilated position. 
If $n=2$ and $K$ and $C$ are in dilated position, then
 Xi, Leng \cite{XiL16} proved 
\eqref{logMconjeqVCK} including the characterization of equality. Actually, \cite{XiL16} even verified Dar's conjecture \eqref{Dar} for convex planar bodies in dilated position (no need for translation in this case).

Let us discuss the existence of the solution of the logaritmic Minkowski problem \eqref{MongeVK} or \eqref{AlexandrovVK}.
Following partial and related results by Andrews \cite{And99}, Chou, Wang \cite{ChW06},
He, Leng, Li \cite{HLL06},
Henk, Sch\"urman, Wills \cite{HSW06}, Stancu \cite{Stancu},
Xiong \cite{Xio10},
the paper B\"or\"oczky, Lutwak, Yang, Zhang \cite{BLYZ13} characterized even cone volume measures
by the so-called subspace concentration condition (i) and (ii) in Theorem~\ref{VKevenchar}.

\begin{theo}
\label{VKevenchar}
There exists an origin symmetric convex body $K\in\mathcal{K}_e^n$ with $\mu=V_K$ for a non-trival finite even Borel measure $\mu$ on $S^{n-1}$ if and only if
\begin{description}
\item[(i)] $\mu(L\cap S^{n-1})\leq \frac{{\rm dim}\,L}{n}\cdot\mu(S^{n-1})$ for any proper linear subspace $L\subset\R^n$;
\item[(ii)] $\mu(L\cap S^{n-1})=\frac{{\rm dim}\,L}{n}\cdot\mu(S^{n-1})$ in (i)
 is equivalent with the existence of a complementary linear subspace $L'\subset\R^n$ with
${\rm supp}\,\mu\subset L\cup L'$.
\end{description}
\end{theo}

We observe that $V_K$ satisfies (ii) if and only if $K=C+C'$ where $C\subset L^\bot$ and  $C'\subset L'^\bot$ compact convex sets. A finite Borel measure $\mu$ on $S^{n-1}$ satisfies the strict subspace concentration condition
if $\mu(L\cap S^{n-1})< \frac{{\rm dim}\,L}{n}\cdot\mu(S^{n-1})$ for any proper linear subspace $L\subset\R^n$.

Given  a non-trival finite even Borel measure $\mu$ on $S^{n-1}$ that is invariant under
$n$ reflections $\Phi_1,\ldots,\Phi_n$ through $n$ independent linear hyperplanes, B\"or\"oczky, Kalantzopoulos \cite{BoK22} proved that $\mu=V_K$ for a convex body $K$ in $\R^n$ invariant under $\Phi_1,\ldots,\Phi_n$ if and only if $\mu$ satisfies  
the subspace concentration condition (i) and (ii) for any proper linear subspace $L\subset\R^n$ invariant under $\Phi_1,\ldots,\Phi_n$. Actually, the statement also holds if $\Phi_1,\ldots,\Phi_n$ are linear reflections (see \cite{BoK22} for details).
For a centered convex body $K\in\mathcal{K}^n$, B\"or\"oczky, Henk \cite{BoH16} (see Henk and Linke \cite{HeL14} for the case of polytopes) verified that $V_K$ satisfies the subspace concentration condition (i) and (ii), but $V_K$ satifies some additional conditions, as well.
 On the other hand, if 
$V_K(L\cap S^{n-1})\geq (1-\varepsilon)\cdot\frac{{\rm dim}\,L}{n}\cdot V(K)$ holds for $K\in\mathcal{K}^n$, a proper linear subspace $L\subset\R^n$
and a small $\varepsilon>0$, then $K$ is close to the sum of two complementary lower dimensional compact convex sets according to B\"or\"oczky, Henk \cite{BoH17}. We note that Freyer, Henk, Kipp \cite{FHK} even verified certain so-called 
Affine Subspace Concentration Conditions for the cone volume measure of centered polytopes.

Much less is known, not even a conjecture about the characteristic properties of a cone volume measure on $S^{n-1}$, not even in the plane. 
Chen, Li, Zhu \cite{CLZ19} proved that if a non-trival finite Borel measure $\mu$ on $S^{n-1}$  satisfies the 
subspace concentration condition (i) and (ii), then $\mu$ is a cone volume measure. On the other hand,
B\"or\"oczky, Heged\H{u}s \cite{BoH15} characterized the restriction of a cone volume measure to a pair of antipodal points.

As Lutwak, Yang, Zhang \cite{LYZ05} conjectured, a Borel probability measure $\mu$ on $S^{n-1}$  satisfies the subspace concentration condition (i) and (ii) if and only if there exists an isotropic linear image $\Phi_*\mu$ for a $\Phi\in{\rm GL}(n)$
according to B\"or\"oczky, Lutwak, Yang, Zhang \cite{BLYZ15} (extending the work by Carlen, Cordero-Erausquin \cite{CaC09} in the discrete case and Klartag  \cite{Kla10}
in the strict subspace concentration condition case). Here the probability measure $\mu$ on $S^{n-1}$ is isotropic if
$n\int_{S^{n-1}}u\otimes u\,d\mu(u)={\rm Id}_n$; or in other words, $\|x\|=n\int_{S^{n-1}}\langle x,u\rangle^2\,d\mu(u)$
for any $x\in\R^n$.\\

Next we turn to the logarithmic Brunn-Minkowski conjecture/inequality. For $\lambda\in(0,1)$ and $K,C\in\mathcal{K}_o^n$, 
we define their logarithmic or $L_0$ linear combination by the formula
$$
(1-\lambda)K+_0\lambda C=
\{x\in\R^n:\langle x,u\rangle\leq h_K(u)^{1-\lambda}h_C(u)^{\lambda}\;\forall u\in S^{n-1}\}.
$$
The $L_0$ linear combination is linear invariant; namely,
if $\Phi\in{\rm GL}(n)$, then $\Phi\Big((1-\lambda)K+_0\lambda C\Big)=(1-\lambda)\Phi(K)+_0\lambda \Phi(C)$. Moreover,
the $L_0$ linear combination is a convex body if $\{h_K=0\}=\{h_C=0\}$ (for example, when $K,C\in\mathcal{K}_{(o)}^n$).
We note that $(1-\lambda)(\alpha K)+_0\lambda (\beta C)=\alpha^{1-\lambda}\beta^\lambda\Big((1-\lambda)K+_0\lambda C\Big)$ for $\alpha,\beta>0$.
The $L_0$ linear combination of polytopes is always a polytope, but the boundary of the $L_0$ linear combination of convex bodies with $C^2_+$ boundaries may contain segments, and hence may not be $C^2_+$. A functional analogue of the $L_0$-addition is presented by Crasta, Fragal\`a \cite{CrF23}.

We observe that
$(1-\lambda)K+_0\lambda C\subset (1-\lambda)K+\lambda C$ for any convex bodies $K$ and $C$ containing the origin interior, but 
  $(1-\lambda)K+_0\lambda C$ might be much smaller than $(1-\lambda)K+\lambda C$. For example, if $a>0$ is large,
$n=2$, $K=[\frac{-1}a,\frac{1}a]\times [-a,a]$ and $C=[-a,a]\times [\frac{-1}a,\frac{1}a]$, then
\begin{equation}
\label{L0example}
\begin{array}{rcl}
\frac12\,K+_0 \frac12\,C&=&[-1,1]^2\\
\frac12\,K+ \frac12\,C&=&\left[-\frac12(a+\frac1a),\frac12(a+\frac1a)\right]^2.
\end{array}
\end{equation}

B\"or\"oczky, Lutwak, Yang, Zhang \cite{BLYZ12} conjectured the following for origin symetric convex bodies, and Martin Henk proposed the version with centered convex bodies (see also \cite{BoK22}).

\begin{conj}[Log-Brunn-Minkowski conjecture]
\label{logBMconj}
If $\lambda\in(0,1)$ and $K$ and $C$ are centered convex bodies in $\R^n$, then
\begin{equation}
\label{logBMconjeq}
V((1-\lambda)K+_0\lambda C)\geq V(K)^{1-\lambda} V(C)^\lambda
\end{equation}
with equality if and only if $K=K_1+\ldots + K_m$ and $C=C_1+\ldots + C_m$ for compact convex sets
	$K_1,\ldots, K_m,C_1,\ldots,C_m$ of dimension at least one where $\sum_{i=1}^m{\rm dim}\,K_i=n$
	and $K_i$ and $C_i$ are dilates, $i=1,\ldots,m$.
\end{conj}

The Log-Brunn-Minkowski Conjecture~\ref{logBMconj} is a significant strengthening of the Brunn-Minkowski inequality for centered convex bodies
(see \eqref{L0example}). Given $K,C\in\mathcal{K}_o^n$, if \eqref{logBMconjeq} holds for all $\lambda\in(0,1)$,
then the Logarithmic Minkowski inequality \eqref{logMconjeq} follows by considering
$\frac{d}{d\lambda}V((1-\lambda)K+_0\lambda C)|_{\lambda=0^+}$ and using Alexandrov's Lemma~\ref{AlexandrovLemma} according to \cite{BLYZ12}. On the other hand, the argument in \cite{BLYZ12} shows that if ${\cal F}$ is any family of convex bodies closed under $L_0$ linear combination, then the Logarithmic Minkowski inequality \eqref{logMconjeq} for all $K,C\in{\cal F}$ is equivalent to the Logarithmic Brunn-Minkowski inequality \eqref{logBMconjeq}
for all $K,C\in{\cal F}$ and $\lambda\in(0,1)$. In particular, the equivalence holds for the family of origin symmetric convex bodies. According to Kolesnikov, Milman \cite{KoM22} and Putterman \cite{Put21}, taking the second derivative of 
$\lambda\mapsto V((1-\lambda)K+_0\lambda C)$ 
for origin symmetric convex bodies $K$ and $C$ in $\R^n$ leads to the conjectured inequality
\begin{equation}
\label{logBMPuttermann}
\frac{V(K,C;1)^2}{V(K)}\geq \frac{n-1}n\,V(K,C;2)+\frac1{n}\int_{S^{n-1}}\frac{h_C^2}{h_K^2}\,dV_K
\end{equation}
that is a strengthened from of Minkowski's second inequality \eqref{Minkowski2ineq}, and is equivalent to the
Log-Brunn-Minkowski conjecture without the charactherization of equality. More precisely,  \cite{KoM22} proves 
that  for a fixed $K\in\mathcal{K}_e^n$  with $C^2_+$ boundary, \eqref{logBMPuttermann} for all smooth $C\in\mathcal{K}_e^n$ is equivalent to a local form of the Log-Brunn-Minkowski around $K$, and \cite{Put21} verifies the global statement.
Therefore, we have the following three equivalent forms of the Log-Brunn-Minkowski conjecture for origin symmetric convex bodies $K$ and $C$ in $\R^n$ (without the charactherization of equality in the case of the third formulation):
\begin{itemize}
\item $V((1-\lambda)K+_0\lambda C)\geq V(K)^{1-\lambda} V(C)^\lambda$ as in \eqref{logBMconjeq};
\item $\int_{S^{n-1}}\log \frac{h_C}{h_K}\,dV_K\geq \frac{V(K)}n\log\frac{V(C)}{V(K)}$ as in \eqref{logMconjeq};
\item $\frac{V(K,C;1)^2}{V(K)}\geq \frac{n-1}n\,V(K,C;2)+\frac1{n}\int_{S^{n-1}}\frac{h_C^2}{h_K^2}\,dV_K$
as in \eqref{logBMPuttermann}.
\end{itemize}

Another equivalent formulation using the Hilbert-Brunn-Minkowski operator \eqref{Hilbert-operator}
 is due to Kolesnikov, Milman \cite{KoM22}, and is discussed at the end of Section~\ref{secLpMinkowski} (see \eqref{lambda1elog}). In addition, Saroglou  \cite{Sar15} verified that the Log-Brunn-Minkowski inequality for any origin symmetric convex bodies is equivalent with the so-called $B$-property: For any origin symmetric convex body $K$ in $\R^n$ and 
$n\times n$ positive definite diagonal matrix $\Phi$, the function $s\mapsto V([-1,1]^n\cap \Phi^sK)$ of $s\in\R$
is log-concave. Yet another equivalent formulation of the Log-Brunn-Minkowski conjecture for origin symmetric convex bodies in $\R^n$
using the "strong $B$-property" is  due to Nayar, Tkocz \cite{NaT20}: For any $N> n$ and $n$-dimensional linear subspace $L$ of $\R^N$, the $n$-volume of $L\cap\prod_{i=1}^N[-e^{t_i},e^{t_i}]$ is a log-concave function of
$(t_1,\ldots,t_N)\in\R^N$. Actually, \cite{NaT20} proves an analogous property of the crosspolytopes.

Saroglou \cite{Sar16} proved that if the log-Brunn-Minkowski Conjecture \eqref{logBMconjeq} holds for any origin symmetric convex bodies $K$ and $C$ and $\lambda\in(0,1)$, then it holds for any even log-concave measure $\mu$ on $\R^n$; namely,
\begin{equation}
\label{logBMlogconv}
\mu((1-\lambda)K+_0\lambda C)\geq \mu(K)^{1-\lambda}\mu(C)^\lambda.
\end{equation}
In turn, the argument in \cite{Sar16} shows that \eqref{logBMlogconv} for the Gaussian meaure $\mu=\gamma_n$ implies
the log-Brunn-Minkowski Conjecture \eqref{logBMconjeq} for origin symmetric convex bodies. Finally, Kolesnikov \cite{Kol20} provides another equivalent formulation of the  log-Brunn-Minkowski Conjecture for origin symmetric convex bodies in terms of
displacement convexity of certain functional of probability measures on the sphere in optimal transportation.

The Log-Brunn-Minkowski Conjecture~\ref{logBMconj} is still open but has been verified in various cases.
In $\R^2$, Conjecture~\ref{logBMconj} is verified by B\"or\"oczky, Lutwak, Yang, Zhang \cite{BLYZ12} for origin symmetric convex bodies, but it is still open for general centered planar convex bodies.
For unconditional convex bodies, the $L_0$ linear combination contains the so-called coordinatewise product (see Saroglou \cite{Sar15}); therefore, the corresponding inequality for the coordinatewise product by
Uhrin \cite{Uhr94}, Bollob\'as, Leader \cite{BoL95} and Cordero-Erausquin, Fradelizi, Maurey \cite{CEFM04},
following from the Pr\'ekopa-Leindler inequality Theorem~\ref{PLn} yields Conjecture~\ref{logBMconj}. 
The  equality case of the Log-Brunn-Minkowski inequality for unconditional convex bodies was clarified by Saroglou \cite{Sar15}
(see also \cite{BoK22}). The Log-Brunn-Minkowski Conjecture~\ref{logBMconj} for convex bodies invariant under reflections through $n$ independent linear hyperplanes is due to B\"or\"oczky, Kalantzopoulos \cite{BoK22}. In addition,
Conjecture~\ref{logBMconj} is proved for complex bodies by Rotem \cite{Rotem}.

Conjecture~\ref{logBMconj} holds for origin symmetric convex bodies in a neighbourhood of a fixed centered ellipsoid $E$; more precisely, for origin symmetric $K$ and $C$  provided $E\subset K,C\subset (1+c_n)E$ where $c_n>0$ depends only on $n$.
In this form, the statement is due to Chen, Huang, Li, Liu \cite{CHLL20} extending the local estimate by 
Kolesnikov, Milman \cite{KoM22}
(an analogues result holds for linear images of $l_q$ balls for $q>2$ if the dimension $n$ is high enough according to \cite{KoM22} and the method of \cite{CHLL20}).
We note that the case when $K$ and $C$ are in a $C^2$ neighbourhood of $E$ was handled earlier by
Colesanti, Livshyts, Marsiglietti \cite{CLM17}.

In some cases when uniqueness of the solution of the Log-Minkowski problem is known, even the stability of the solution has been established. For example, B\"or\"oczky, De \cite{BoD21b} established this among convex bodies invariant under $n$ given reflections through linear hyperplanes. Concerning Firey's classical result that the only origin symmetric solution of the Log-Minkowski problem \eqref{MongeVK}  with constant $f$ is the centered ball, Ivaki \cite{Iva22} verified a stability version.
Next B\"or\"oczky, Saroglou \cite{BoS} proved the uniqueness of the solution if the possibly non-even $f$ in \eqref{MongeVK} is $C^\alpha$ close to a constant function (the case $n=3$ was handled earlier by Chen, Feng, Liu \cite{CFL22}). 

If $n=2$ and $K$ and $C$ are in dilated position \eqref{dilated}, then
 Xi, Leng \cite{XiL16} proved 
\eqref{logBMconjeq} for any $\lambda\in(0,1)$ including the characterization of equality using the same method as in the case of
the planar Dar conjecture. It is an intriguing question what the relation between Dar's conjecture \eqref{Dar} and the Log-Brunn-Minkowski Conjecture~\ref{logBMconj} is, whether one of them implies the other for origin symmetric convex bodies. 

We note that there exist $\eta_2>\eta_1>0$ depending on $n$ such that if  $\lambda\in(0,1)$ and $K$ and $C$ are centered convex bodies in $\R^n$, then 
\begin{equation}
\label{logBMest}
\eta_1V(K)^{1-\lambda} V(C)^\lambda\leq V((1-\lambda)K+_0\lambda C)\leq \eta_2V(K)^{1-\lambda} V(C)^\lambda,
\end{equation}
which estimates indicate why proving the Log-Brunn-Minkowski Conjecture~\ref{logBMconj} is so notoriously difficult.
Conjecture~\ref{logBMconj} states that $\eta_1=1$, but here we only verify that $\eta_1=n^{-n}$ and 
$\eta_2=n^{3n/2}$ work. According to 
 Kannan, Lov\'asz, Simonovits \cite{KLS95}, there exist centered ellipsoids $E'\subset K$ and $E\subset C$ such that
$K\subset nE'$ and $C\subset nE$. After a linear transform, we may assume that $E'=B^n_2$ and $E$ is unconditional. 
Since Conjecture~\ref{logBMconj} holds for the unconditional convex bodies $B^n_2$ and $E$, we deduce that 
$\eta_1=n^{-n}$ works in \eqref{logBMest}.
For the upper bound, let $a_1,\ldots,a_n$ be the half axes of $E$, and hence 
$C\subset \widetilde{C}=\prod_{i=1}^n[-na_i,na_i]$ and $K\subset \widetilde{K}=[-n,n]^n$ with
$V(\widetilde{C})\leq n^{3n/2}V(C)$ and $V(\widetilde{K})\leq n^{3n/2}V(K)$. Since
$V((1-\lambda)\widetilde{K}+_0\lambda \widetilde{C})= V(\widetilde{K})^{1-\lambda} V(\widetilde{C})^\lambda$,
it follows that $\eta_2=n^{3n/2}$ works in \eqref{logBMest}.

The validity of the Log-Minkowski (or Log-Brunn-Minkowski) Conjecture is also supported by the fact that various consequences of it has been verified. For example, the $L_p$-Minkowski Conjecture has been proved when $p\in(0,1)$ is close to $1$ (see Theorem~\ref{pcloseto1}). Nest we turn to results about  the canonical Gaussian probability measure $\gamma_n$ on $\R^n$. One possible consequence of the Log-Brunn-Minkowski Conjecture~\ref{logBMconj}
is the earlier celebrated "$B$-inequality" by Cordero-Erausquin, Fradelizi, Maurey \cite{CEFM04} stating that $\gamma_n(e^tK)$ is a log-concave function of $t\in\R$ for any origin symmetric $K\in\mathcal{K}^n$.
Next, the Gardner-Zvavitch conjecture in \cite{GaZ10} stated that
if $K$ and $C$ are origin symmetric convex bodies in $\R^n$, then
\begin{equation}
\label{BM-Gaussian-measure}
\gamma_n((1-\lambda)K+\lambda C)^{\frac1n}\geq (1-\lambda)\gamma_n(K)^{\frac1n}+\lambda\gamma_n(C)^{\frac1n}.
\end{equation}
It was proved by Livshyts, Marsiglietti, Nayar, Zvavitch \cite{CoL20},
that the log-Brunn-Minkowski conjecture would imply the Gardner-Zvavitch conjecture.
After various attempts, the conjecture was finally verified by 
Eskenazis, Moschidis \cite{EsM21} not much before that,
Kolesnikov, Livshyts \cite{KoL21} verified that if the exponents $\frac1n$ in \eqref{BM-Gaussian-measure} are changed into 
$\frac1{2n}$, then this modified Gardner-Zvavitch conjecture holds for any pair of centered convex bodies $K$ and $C$.

We note that independently of the log-Brunn-Minkowski conjecture, various Brunn-Minkowski type inequalities have been proved and conjectured for the Gaussian measure, the most famous ones being the Ehrhardt inequality and the Gaussian isoperimetric inequality (see Livshyts \cite{Liva}).

Colesanti, Livshyts, Marsiglietti \cite{CLM17} conjectured the following generalization of the Gardner-Zvavitch conjecture:
 If $\mu$ is an even log-concave measure on $\R^n$, then
\begin{equation}
\label{BMlog-concave-measure}
\mu((1-\lambda)K+\lambda C)^{\frac1n}\geq (1-\lambda)\mu(K)^{\frac1n}+\lambda\mu(C)^{\frac1n}
\end{equation}
holds for any origin symmetric convex bodies $K$ and $C$.
According to  Livshyts, Marsiglietti, Nayar, Zvavitch \cite{CoL20},
 the Log-Brunn-Minkowski Conjecture~\ref{logBMconj} would imply the conjecture \eqref{BMlog-concave-measure}.
Cordero-Erausquin, Rotem \cite{CER} proved \eqref{BMlog-concave-measure} if $\mu$ is a rotationally symmetric log-concave measure.
In addition, 
Livshyts \cite{Livb} verified that \eqref{BMlog-concave-measure} holds for
any even log-concave measure on $\R^n$ and origin symmetric convex bodies $K$ and $C$
if the exponents $\frac1n$ in \eqref{BMlog-concave-measure} are changed into 
$n^{-4-o(1)}$.

\section{Lutwak's $L_p$-Minkowski theory}
\label{secLpMinkowski}

The rapidly developing new  $L_p$-Brunn-Minkowski theory (where $p=1$ is the classical case and $p=0$ corresponds to the cone-volume measure)
initiated by Lutwak \cite{Lut93,Lut93a,Lut96}, has become main research area in modern convex geometry and geometric analysis. For $p\in \R$ and $K\in\mathcal{K}_o^n$, the $L_p$-surface area measure $S_{K,p}$ on $S^{n-1}$ is defined by
\begin{equation}
\label{SKp}
dS_{K,p}=h_K^{1-p}\,dS_K
\end{equation}
where if $p>1$ and $o\in\partial K$, then we assume that $S_K(\{h_K=0\})=0$. In particular, $S_{K,1}=S_K$ and $S_{K,0}=nV_K$. For $p\in \R$, the Monge-Amp\'ere equation on $S^{n-1}$ corresponding to the $L_p$-Minkowski problem is
\begin{equation}
\label{MongeSKp}
\begin{array}{rcl}
\det(\nabla^2 h+h\,{\rm Id})= h^{p-1}f&\mbox{ if }&p>1\\[1ex]
h^{1-p}\det(\nabla^2 h+h\,{\rm Id})= f&\mbox{ if }&p\leq 1
\end{array}
\end{equation}
where $f\in L_1(S^{n-1})$ is non-negative with $\int_{S^{n-1}}fd\HH>0$,
and for a finite non-trivial Borel measure $\mu$ on $S^{n-1}$, a convex body $K\in\mathcal{K}_o^n$ is an Alexandrov solution
of the $L_p$-Minkowski problem if
\begin{equation}
\label{AlexandrovSKp}
\begin{array}{rcl}
dS_K= h_K^{p-1}\,d\mu&\mbox{ if }&p>1\\[1ex]
h_K^{1-p}dS_K= d\mu&\mbox{ if }&p\leq 1.
\end{array}
\end{equation}

If $p>1$ and $p\neq n$, then Hug, Lutwak, Yang, Zhang \cite{HLYZ2} (improving on Chou, Wang \cite{ChW06} ) prove that
\eqref{AlexandrovSKp} has an Alexandrov solution if and only if the $\mu$ is not concentrated onto any closed hemisphere,
and the solution is unique. If in addition $p>n$, then the unique solution of \eqref{AlexandrovSKp} satifies 
$o\in{\rm int}\,K$, and hence $S_{K,p}=\mu$.
However, examples in \cite{HLYZ2} show that if $1<p<n$, then it may happen that the density function $f$ is a positive continuous in \eqref{MongeSKp} and $o\in\partial K$ holds for the unique  Alexandrov solution.
If $p=n$, then $S_{K,n}=S_{\lambda K,n}$ holds for $\lambda>0$; therefore, all what is known (see \cite{HLYZ2}) is that for any
measure $\mu$ not concentrated onto any closed hemisphere, there exists a convex body $K\in\mathcal{K}_o^n$ and $c>0$ such that $\mu=c\cdot S_{K,n}$.

The case $p=1$ is the classical Minkowski problem (see Section~\ref{secMinkowski}), and the case 
case $p=0$ is the logarithmic Minkowski problem (see Section~\ref{secLogMinkowski}). 

If $p\in(0,1)$ and the
measure $\mu$ is not concentrated onto any great subsphere, then Chen, Li, Zhu \cite{CLZ17} prove that there 
exists an Alexandrov solution $K\in\mathcal{K}_o^n$ of \eqref{AlexandrovSKp} with $S_{K,p}=\mu$. For $p\in(0,1)$, complete characterization of $L_p$ surface area measures is only known
if $n=2$ by B\"or\"oczky, Trinh \cite{BoT17}; namely, a finite non-trivial Borel measure $\mu$ on $S^1$ is
an $L_p$ surface area measure if and only if ${\rm supp}\,\mu$ does not consists of a pair of antipodal points. Finally, let 
$p\in(0,1)$ and $n\geq 3$, and let us assume that $1\leq{\rm dim}\,L\leq n-1$ where $L$ is the linear hull of ${\rm supp}\,\mu$ in $\R^n$.
If ${\rm supp}\,\mu$ is contained in a closed hemisphere centered at a point of $L\cap S^{n-1}$, then $\mu$ is an
$L_p$ surface area measure according to Bianchi, B\"or\"oczky, Colesanti, Yang \cite{BBCY19}. On the other hand, 
Saroglou \cite{Sar} proved that if
$\mu(\omega)$ is the Lebesgue measure of $\omega\cap L$ for any Borel $\omega\subset S^{n-1}$, then
$\mu$ is not a $L_p$ surface area measure.

If $-n<p<0$ and $f\in L_{\frac{n}{n+p}}(S^{n-1})$ in \eqref{MongeSKp}, then \eqref{MongeSKp} has a solution according to
Bianchi, B\"or\"oczky, Colesanti, Yang \cite{BBCY19}. If $p<0$ and the $\mu$ in \eqref{AlexandrovSKp} is discrete satisfying that $\mu$ is not concentrated on any closed hemisphere and any $n$ unit vectors in the support of $\mu$ are independent, then
\cite{Zhu17} manages to solve the $L_p$-Minkowski problem.

The $p=-n$ case of the $L_p$-Minkowski problem is the critical case because its link with  the ${\rm SL}(n)$ invariant
centro-affine curvature.
If $K\in \mathcal{K}_{(o)}^n$ has $C^2_+$ boundary, then its centro-affine curvature at $u\in S^{n-1}$ is
$\kappa_0(K,u)=\frac{\kappa(K,u)}{h_K(u)^{n+1}}$ where $\kappa(K,u)$ is the Gaussian curvature at the point with exterior normal $u$.
It is well known to be ${\rm SL}(n)$ invariant  in the sense that
$\kappa_0(\Phi K,u)=\kappa_0(K,\frac{\Phi^{t} u}{\|\Phi^{t} u\|})$ for $\Phi\in{\rm SL}(n)$
 (see Hug \cite{Hug96} or Ludwig \cite{Lud10}). It follows from \eqref{SKGauss} that if $K\in \mathcal{K}_{(o)}^n$ has $C^2_+$ boundary, then $dS_{K,-n}(u)=\kappa_0(K,u)^{-1}\,d\HH(u)$; therefore, solving the $L_p$-Minkowski problem
\eqref{MongeSKp} for $p=-n$ and positive $C^\alpha$ function $f$ is equivalent to reconstructing a convex body 
$K\in \mathcal{K}_{(o)}^n$ from its centro-affine curvature function.

All in all, the centro-affine ($L_{-n}$) Minkowki problem is wide open. If $p=-n$ and the $f$ in \eqref{MongeSKp} is unconditional and satisfies certain additional technical conditions, then
Jian, Lu, Zhu \cite{JLZ16} verify the existence of a solution of  \eqref{MongeSKp}.
Moreover the paper Li, Guang, Wang \cite{LGWc} solves a variant of the centro-affine Minkowki problem. On the other hand, Chou, Wang \cite{ChW06} prove an implicit condition on possible functions $f$ in \eqref{MongeSKp} such that in $f^{-1}$ is a centro affine curvature (see also \cite{BBCY19}), and
Du \cite{Du21} construct an explicit example of a positive $C^\alpha$ $f$ such that \eqref{MongeSKp} has no solution when $p=-n$. 

In the super-critical case $p<-n$, Li, Guang, Wang \cite{LGWa} have recently achieved a breakthrough by proving that 
for any positive $C^2$ $f$, there exists a $C^4$ solution of  \eqref{MongeSKp}. 
In addition, 
\cite{LGWa} verify that if $p<-n$ and $1/c<f<c$ for a constant $c>1$ in \eqref{MongeSKp}, then there exists a $C^{1,\alpha}$ Alexandrov solution $h_K|_{S^{n-1}}$ satisfying \eqref{AlexandrovSKp} where $d\mu=f\,d\HH$.
In their paper,  Guang, Li, Wang \cite{LGWa} combine a flow argument with homology calculations.
On the other hand, Du \cite{Du21} construct a non-negative $C^\alpha$ function $f$ that is positive everywhere but a fixed pair of antipodal points and \eqref{MongeSKp} has no solution, not even in Alexandrov sense.
It is not surprising that the flow argument works in the super-critical case, as Milman \cite{Mila} points out the limitations of the variational argument in this case.
For a discrete measure $\mu$ satisfying that $\mu$ is not concentrated on any closed hemisphere and any $n$ unit vectors in the support of $\mu$ are independent, Zhu \cite{Zhu17} solves the $L_p$-Minkowski problem \eqref{AlexandrovSKp} for $p<0$.

If $p>-n$, then while flow arguments are also known (see e.g. Bryan, Ivaki, Scheuer \cite{BIJ19}), the most common argument to find a solution of  \eqref{MongeSKp} is based on the variational method; namely, one considers the infimum of $\int_{S^{n-1}}h_C^pf\,d\HH$ for a suitable family of convex bodies $C\in\mathcal{K}_{(0)}^n$ with $V(C)=1$ when $f$ is positive and continuous (see \cite{BBCY19} or Chou, Wang \cite{ChW06}). The existence of some minimizer $C_0$ follows via the Blaschke-Santal\'o inequality \eqref{Blaschke-Santalo} as $p>-n$, and the fact that $dS_{C_0,p}=\lambda f\,d\HH$ for some constant factor $\lambda>0$ follows via the Alexandrov Lemma~\ref{AlexandrovLemma}. The case of more general measures than the ones with positive continuous density functions follows by approximation. For the variational approach, it is also common to use discrete measures on $S^{n-1}$ (corresponding to polytopes, see \cite{HLYZ2,Zhu15,BHZ16,Zhu17}). 

Concerning the smoothness of the solution of the $L_p$-Minkowski problem \eqref{MongeSKp}, if $f$ is positive and $C^\alpha$ and $h$ is positive (equivalently, $o\in{\rm int}\,K$ for the corresponding convex body $K$), then $h$ is $C^{2,\alpha}$ by Cafarelli \cite{Caf90a,Caf90b}
(see \cite{BoF19}, \cite{BBC20}). Assuming that $f$ is positive and continuous, it is known that $o\in{\rm int}\,K$
if $p\geq n$ (see \cite{HLYZ2}) or if $p\leq 2-n$ (see \cite{BBC20}). On the other hand, if $2-n<p<n$, $p\neq 1$, then
there there exists positive $C^\alpha$ function $f$ on $S^{n-1}$ such that $o\in\partial K$ holds for the
Alexandrov solution of \eqref{AlexandrovSKp} with $d\mu=f\,d\HH$, see \cite{HLYZ2}
if $1<p<n$ and \cite{BBC20} if $2-n<p<1$. Additional results about the smoothness of the solution 
are provided by \cite{BBC20} in the case $2-n<p<1$.

Now we discuss the uniqueness of the solution of the $L_p$-Minkowski problem \eqref{AlexandrovSKp}.
As we have seen, if $p>1$ and $p\neq n$, then Hug, Lutwak, Yang, Zhang \cite{HLYZ2} proved that the Alexandrov solution of 
the $L_p$-Minkowski problem \eqref{AlexandrovSKp} is unique. However, if $p<1$, then
the solution of the $L_p$-Minkowski problem \eqref{MongeSKp} may not be unique even if $f$ is positive and continuous.
Examples are provided by Chen, Li, Zhu \cite{CLZ17,CLZ19} if $p\in[0,1)$, and Milman \cite{Mila} shows that for any $C\in\mathcal{K}_{(0)}$, one finds $q\in [-n,1)$ where $q=-n$ exactly when $C$ is a centered ellipsoid such that if $p<q$, then there exist multiple solutions of
the $L_p$-Minkowski problem \eqref{AlexandrovSKp} with $\mu=S_{C,p}$; or in other words, there exists 
$K\in\mathcal{K}^n_{(0)}$ with $K\neq C$ and $S_{K,p}=S_{C,p}$. In addition, 
Jian, Lu, Wang \cite{JLW15} and Li, Liu, Lu \cite{LLL22} prove that for any $p<0$, there exists positive even $C^\infty$ function $f$ with rotational symmetry
such that the $L_p$-Minkowski problem \eqref{MongeSKp} has multiple positive  even $C^\infty$ solutions.
We note that in the case of the centro-affine Minkowski problem $p=-n$, Li \cite{Li19} even verified the possibility of existence of infinitely many solutions without affine equivalence, and Stancu \cite{Sta22} proved that if an origin symmetric convex body $K$ with $C^\infty_+$ boundary is a unique solution to the $L_p$-Minkowski problem \eqref{AlexandrovSKp} up to linear equivalence for $p=-n$
with $\mu=S_{K,-n}$, then it is  a unique solution for $p=0$
with $\mu=S_{K,0}$.

The case when $f$ is a constant function in the $L_p$-Minkowski problem \eqref{MongeSKp} has received a special attention
since Firey \cite{Fir74}. Through the work of Lutwak \cite{Lut93a}, Andrews \cite{And99}, Andrews, Guan, Ni \cite{AGN16} and Brendle, Choi, Daskalopoulos \cite{BCD17}, it has been clarified that the only solutions are centered balls if $p>-n$, centered ellipsoids if $p=-n$, and there are several solutions if $p<-n$.
See Crasta, Fragal\'a \cite{CrF23}, Ivaki, Milman \cite{IvM} and Saroglou \cite{Sar22} for novel approaches.
 Stability versions of these results have been obtained by Ivaki \cite{Iva22}, sometimes in the even case for certain ranges of $p$.
Uniqueness of the solution of \eqref{MongeSKp} if $f$ is $C^{\alpha}$ close to $1$ (without the evenness assumption) is proved if $p\in[0,1)$ by B\"or\"oczky, Saroglou \cite{BoS}.

In particular, concerning uniqueness, the a major question left open is the
uniqueness of even solutions of the $L_p$-Minkowski problem \eqref{MongeSKp} when 
$f$ is an positive even $C^\infty$ function and $p\in[0,1)$. In the case of $p=0$, this is Lutwak's
Log-Minkowski Conjecture~\ref{logMinkowskiUniqCinfty}. If $p\in(0,1)$, it is also conjectured
that  the $L_p$-Minkowski problem \eqref{MongeSKp} has a unique even solution for any positive, $C^\infty$ and even $f$.
More generally, we have the following conjecture (see B\"or\"oczky, Lutwak, Yang, Zhang \cite{BLYZ12} for origin symmetric bodies).

\begin{conj}[$L_p$-Minkowski Conjecture \#1]
\label{LpMinkowskiUniqCinfty}
If $p\in(0,1)$ and $K$ and $C$ are centered convex bodies in $\R^n$ with $S_{K,p}=S_{C,p}$, then $K=C$.
\end{conj}

Before presenting what is known about the $L_p$-Minkowski conjecture, let us discuss its relation to the 
$L_p$-Brunn-Minkowski theory for $p\geq 0$. More precisely, the cases $p=0$ and $p=1$ have been discussed in 
Sections~\ref{secMinkowski} and \ref{secLogMinkowski}.

For $p>0$, $\alpha,\beta>0$ and $K,C\in\mathcal{K}_o^n$, 
we define the $L_p$ linear combination by the formula
$$
(1-\lambda)K+_p\lambda C=
\{x\in\R^n:\langle x,u\rangle^p\leq \alpha\,h_K(u)^p+\beta\,h_C(u)^p\;\forall u\in S^{n-1}\}.
$$
The $L_p$ linear combination is linear invariant; namely,
if $\Phi\in{\rm GL}(n)$, then $\Phi\Big(\alpha\,K+_p\beta C\Big)=\alpha\,\Phi(K)+_p\beta\, \Phi(C)$. 
If $p\in (0,1)$, then the $L_p$ linear combination of polytopes is always a polytope, but the boundary of the $L_p$ linear combination of convex bodies with $C^2_+$ boundaries may contain segments, and hence may not be $C^2_+$.
On the other hand, if $p>1$, then for any $\alpha,\beta>0$, Minkowski's inequality yields that
$h_{\alpha\,K+_p\beta\,C}^p=\alpha\,h_{K}^p+\beta\,h_{C}^p$, as the $L_p$ linear combination was defined by Firey \cite{Fir62} in this case. According to Firey \cite{Fir62}, if $p>1$ and $K,C\in\mathcal{K}_o^n$, then the Brunn-Minkowski inequality yields the 
$L_p$-Brunn-Minkowski inequality
\begin{equation}
\label{LpBMeq1}
V(\alpha\,K+_p\beta\, C)^{\frac{p}n}\geq \alpha\,V(K)^{\frac{p}n} +\beta\,V(C)^{\frac{p}n}
\end{equation}
for any $\alpha,\beta>0$ with equality if and only if $K$ and $C$ are dilates;
or equivalently,
\begin{equation}
\label{LpBMeq2}
V((1-\lambda)K+_p\lambda C)\geq V(K)^{1-\lambda} V(C)^\lambda
\end{equation}
for $\lambda\in(0,1)$ with equality if and only if $K=C$. 

For $p>0$ and $K,C\in\mathcal{K}_{(o)}^n$, analogously to the classical mixed volumes, Lutwak \cite{Lut93} introduced the $L_p$ mixed volume
$$
V_p(K,C)=\frac{p}n\lim_{t\to 0^+}\frac{V(K+_pt\, C)-V(K)}t
=\frac1n\int_{S^{n-1}}h_C^p\,dS_{K,p}=
\int_{S^{n-1}}\frac{h_C^p}{h_K^p}\,dV_K,
$$
and hence $V_1(K,C)=V(K,C;1)$.
Considering the first derivative 
of $\lambda\mapsto V((1-\lambda)K+_p\lambda C)^{\frac{p}n}$ yields the $L_p$-Minkowski inequality
\begin{equation}
\label{LpMinkowskieq}
V_p(K,C)\geq V(K)^{\frac{n-p}n}V(C)^{\frac{p}n}
\end{equation}
for $p>1$ and $K,C\in\mathcal{K}_{(o)}^n$ with equality if and only if $K$ and $C$ are dilates. An equivalent form is that 
\begin{equation}
\label{LpMinkowskieq0}
\int_{S^{n-1}}h_C^p\,dS_{K,p}\geq \int_{S^{n-1}}h_K^p\,dS_{K,p}
\end{equation}
if $p>1$, $K,C\in\mathcal{K}_{(o)}^n$ and $V(K)=V(C)$ with equality if and only if $K=C$. 

We recall that that the Brunn-Minkowski inequality \eqref{BrunnMinkowski} holds for
 bounded Borel subsets $K$ and $C$ of $\R^n$, as well. When $p>1$, the $L_p$-Brunn-Minkowski inequality has been also extended to certain families of non-convex sets by
Zhang \cite{Zha99}, Ludwig, Xiao, Zhang \cite{LXZ11} and Lutwak, Yang, Zhang \cite{LYZ12}.

If $p\in(0,1)$, then translating a cube shows that neither $L_p$-Brunn-Minkowski inequality, nor 
$L_p$-Minkowski inequality hold for general $K,C\in\mathcal{K}_{(o)}^n$. However, 
B\"or\"oczky, Lutwak, Yang, Zhang \cite{BLYZ12} conjecture that they hold for at least origin symmetric convex bodies
(see B\"or\"oczky, Kalantzopoulos \cite{BoK22} for centered convex bodies).

\begin{conj}[$L_p$-Minkowski conjecture \#2]
\label{LpMconj}
If $p\in(0,1)$, then \eqref{LpMinkowskieq}; or equivalently, 
\eqref{LpMinkowskieq0} hold for centered $K,C\in\mathcal{K}^n$.
\end{conj}

\begin{conj}[$L_p$-Brunn-Minkowski conjecture]
\label{LpBMconj}
If $p\in(0,1)$, then \eqref{LpBMeq1}; or equivalently, 
\eqref{LpBMeq2} hold for centered $K,C\in\mathcal{K}^n$.
\end{conj}

The fact that the forms Conjectures~\ref{LpMinkowskiUniqCinfty} and  \ref{LpMconj} (including the characterization of equality)
of the $L_p$-Minkowski conjecture are equivalent follows from \eqref{LpMinkowskieq0} and the variational method as described above.

According to the Jensen inequality, $(1-\lambda)K+_q\lambda \subset (1-\lambda)K+_p\lambda C$ for $p>q\geq 0$.
It follows for example {\it via} \eqref{LpBMeq2} (or {\it via} \eqref{lambda1ep}) that if $0\leq q<p<1$, then the $L_q$-Brunn-Minkowski conjecture (or equivalently the $L_q$-Minkowski conjecture) yields the
$L_p$-Brunn-Minkowski conjecture (or equivalently the $L_q$-Minkowski conjecture). In particular, the 
 $L_p$-Brunn-Minkowski Conjecture~\ref{LpBMconj} for $p\in(0,1)$ is a strengthening of the Brunn-Minkowski inequality for centered convex bodies on the one hand, and follows from the Log-Brunn-Minkowski Conjecture~\ref{logBMconj} on the other hand. In addition Kolesnikov-Milman \cite{KoM22} prove that knowing the $L_p$-Minkowski inequality \eqref{LpMinkowskieq} for some
$p\in(0,1)$ yields even the characterization of the equality  case for the $L_q$-Minkowski inequality when $q\in(p,1)$.

Let $p\in(0,1)$. Given $K,C\in\mathcal{K}_o^n$, if \eqref{LpBMeq1} holds for all $\alpha,\beta>0$,
then the $L_p$-Minkowski inequality \eqref{LpMinkowskieq} follows by considering
$\frac{d}{d\lambda}V((1-\lambda)K+_0\lambda C)|_{\lambda=0^+}$ and using Alexandrov's Lemma~\ref{AlexandrovLemma} according to \cite{BLYZ12}. On the other hand, the argument in \cite{BLYZ12} shows that if ${\cal F}$ is any family of convex bodies closed under $L_p$ linear combination, then the $L_p$-Minkowski inequality \eqref{LpMinkowskieq} for all $K,C\in{\cal F}$ is equivalent to the $L_p$ Brunn-Minkowski inequality \eqref{LpBMeq1}
for all $K,C\in{\cal F}$ and $\alpha,\beta>0$. In particular, the equivalence holds for the family of origin symmetric convex bodies. According to Kolesnikov, Milman \cite{KoM22} and Putterman \cite{Put21}, taking the second derivative of $\lambda\mapsto V((1-\lambda)K+_p\lambda C)^{\frac{p}n}$ 
for origin symmetric convex bodies $K$ and $C$ in $\R^n$ leads to the conjectured inequality
\begin{equation}
\label{LpBMPuttermann}
\frac{V(K,C;1)^2}{V(K)}\geq  \frac{n-1}{n-p}\,V(K,C;2)+ \frac{1-p}{n-p}\,\int_{S^{n-1}}\frac{h_C^2}{h_K^2}\,dV_K
\end{equation}
that is again a strengthened from of Minkowski's second inequality \eqref{Minkowski2ineq},
and is equivalent to the
$L_p$ Brunn-Minkowski conjecture without the charactherization of equality. More precisely,  \cite{KoM22} proves 
that  for a fixed $K\in\mathcal{K}_e^n$  with $C^2_+$ boundary, \eqref{logBMPuttermann} for all smooth 
$C\in\mathcal{K}_e^n$ is equivalent to a local form of the $L_p$ Brunn-Minkowski around $K$, and \cite{Put21} verifies the global statement. We note that van Handel \cite{vHa} presents an approach relating the equality case of 
\eqref{LpBMPuttermann} to the equality case of \eqref{LpMinkowskieq}
for a fixed $K\in\mathcal{K}_e^n$  with $C^2_+$ boundary.

In summary, we have the following three equivalent forms of the $L_p$-Brunn-Minkowski conjecture for 
$p\in(0,1)$ and origin symmetric convex bodies $K$ and $C$ in $\R^n$ (without the charactherization of equality in the case of the third formulation):
\begin{itemize}
\item $V((1-\lambda)K+_p\lambda C)
\geq V(K)^{1-\lambda}V(C)^{\lambda}$ for $\lambda\in(0,1)$;
\item $V_p(K,C)\geq V(K)^{\frac{n-p}n}V(C)^{\frac{p}n}$;
\item $\frac{V(K,C;1)^2}{V(K)}\geq \frac{n-1}{n-p}\,V(K,C;2)+\frac{1-p}{n-p}\int_{S^{n-1}}\frac{h_C^2}{h_K^2}\,dV_K$.
\end{itemize}

Let us dicuss the cases when Conjectures~\ref{LpMinkowskiUniqCinfty}, \ref{LpMconj} and \ref{LpBMconj} have been verified.
They have been verified in the planar $n=2$ case by B\"or\"oczky, Lutwak, Yang, Zhang \cite{BLYZ12}.
 The most spectacular result is due to the combination of the local result by Kolesnikov, Milman \cite{KoM22} and
the local to global approach based on Schrauder estimates in PDE by Chen, Huang, Li, Liu \cite{CHLL20} (see Puttermann \cite{Put21} for an Alexandrov-type argument for the local to global approach) is that
the $L_p$-Minkowski and $L_p$-Brunn-Minkowski conjectures hold for origin symmetric convex bodies if $p\in(0,1)$ is close to $1$.

\begin{theo}
\label{pcloseto1}
If $n\geq 3$ and $p\in(p_n,1)$ where $0<p_n<1-\frac{c}{n(\log n)^{10}}$ for an absolute constant $c>0$, then
the $L_p$-Brunn-Minkowski and $L_p$-Minkowski conjectures 
\eqref{LpBMeq1}, \eqref{LpBMeq2}, \eqref{LpMinkowskieq} and
\eqref{LpMinkowskieq0} hold for $K,C\in\mathcal{K}_e^n$ including the characterization of the equality cases. 
\end{theo}

The paper  Kolesnikov, Milman \cite{KoM22} provides an explicit estimate for $p_n$, depending on the Cheeger or Poincar\'e constants subject of the celebrated Kannan, Lov\'asz, Simonovits conjecture \cite{KLS95}. Our estimate for $p_n$
comes from the upper bound $c(\log n)^5$ by Klartag, Lehec \cite{KlL} for  the Poincar\'e constant where $c>0$ is a absolute constant.

Otherwise, the known cases of the $L_p$-Minkowski and $L_p$-Brunn-Minkowski conjectures for origin symmetric bodies follow from the known cases of the Log-Minkowski and   Log-Brunn-Minkowski conjectures.
Let $p\in (0,1)$ and $n\geq 3$. Then \eqref{LpBMeq1}, \eqref{LpBMeq2}, \eqref{LpMinkowskieq} and
\eqref{LpMinkowskieq0} hold if $K$ and $C$ are invariant under reflections through fixed $n$ independent linear hyperplanes
({\it cf.}  \cite{BoK22}) and if $K$ and $C$ are origin symmetric complex bodies ({\it cf.}  \cite{Rotem}).
In addition, the $L_p$-Minkowski conjecture \eqref{LpMinkowskieq} and
\eqref{LpMinkowskieq0} hold for $K,C\in\mathcal{K}_e^n$ (together with characterization of equality if $\partial K$ is $C^3_+$) if either $K$ is a zonoid according to  \cite{vHa}, 
or there exists a centered ellipsoid $E$ with $E\subset K\subset (1+c_n)E$ where $c_n>0$ depends only on $n$
according to \cite{CHLL20}  (an analogues result holds for linear images of $l_q$ balls for $q>2$ if the dimension $n$ is high enough according to \cite{KoM22}).

Concerning the $L_p$ Brunn-Minkowski conjecture,
Hosle, Kolesnikov, Livshyts \cite{HKL21} and Kolesnikov, Livshyts \cite{KoL22}
 present certain natural generalizations and approaches.

In the final part of Section~\ref{secLpMinkowski}, we discuss how David Hilbert's elegant operator theoretic proof of the Brunn-Minkowski inequality has lead to recent new approaches initiated by Kolesnikov, Milman  \cite{KoM22} towards the 
$L_p$-Minkowski conjecture (see also Putterman \cite{Put21} and van Handel \cite{vHa}).
Here we present Kolesnikov and Milman's version  of
the Hilbert-Brunn-Minkowski operator based on \cite{KoM22} because this modified
operator $\mathcal{L}_K$ intertwines with linear transformations ({\it cf.} Theorem~5.8 in \cite{KoM22}).

The mixed discriminat $D_\ell(B_1,\ldots,B_\ell)$ of $\ell$ positive definite $\ell\times\ell$ matrices can be defined  via the identity
\begin{equation}
\label{Mixeddisc}
\mbox{$\det_\ell$}\left(\lambda_1A_1+\ldots+\lambda_mA_m\right)=\sum_{i_1,\ldots,i_\ell=1}^m D_\ell(A_{i_1},\ldots,A_{i_\ell})
\cdot \lambda_{i_1}\cdot\ldots\cdot\lambda_{i_\ell}
\end{equation}
for $\lambda_1,\ldots,\lambda_m\in\R$ and positive definite $\ell\times\ell$ matrices
$A_1,\ldots,A_m$  where 
$D_\ell(A_{i_1},\ldots,A_{i_\ell})> 0$ is symmetric in its variables and $D_\ell(A,\ldots,A)=\det A$
 (see van Handel, Shenfeld \cite{vHS19,vHS} or Kolesnikov, Milman \cite{KoM22}). The similarity between \eqref{Mixedvol} and \eqref{Mixeddisc} is not a coincidence as
$$
V(K_1,\ldots,K_n)=\frac1n\int_{S^{n-1}} h_{K_n} D_{n-1}\left(D^2h_{K_1},\ldots,D^2h_{K_{n-1}}\right)\,d\HH
$$
for 
$K_1,\ldots,K_n\in\mathcal{K}^n$ with $C^2_+$ boundary. 

For $K\in\mathcal{K}_{(o)}^n$ with $C^2_+$ boundary,  Kolesnikov, Milman \cite{KoM22} defines the
Hilbert-Brunn-Minkowski operator $\mathcal{L}_K:C^2(S^{n-1})\to C^2(S^{n-1})$ by the formula
\begin{equation}
\label{Hilbert-operator}
\mathcal{L}_K f=\frac{D_{n-1}\left(D^2(fh_K),D^2h_{K},\ldots,D^2h_{K}\right)}{D_{n-1}\left(D^2h_K,\ldots,D^2h_{K}\right)}-f.
\end{equation}
Following Hilbert's footsteps, \cite{KoM22} verifies that  the operator $\mathcal{L}_K$ is elliptic,  and hence admits a unique self-adjoint extension
in $L^2(dV_K)$, and has discrete spectrum. The operator $-\mathcal{L}_K$ is positive semi-definite, its smallest eigenvalue is
$\lambda_0(\mathcal{L}_K)=0$ whose eigenspace consists of the constant functions. As Hilbert (see also van Handel, Shenfeld \cite{vHS19,vHS} or Kolesnikov, Milman \cite{KoM22}) proved, the next eigenvalue is 
$\lambda_1(-\mathcal{L}_K)=1$ corresponding to the $n$-dimensional eigenspace spanned by the linear functions; moreover, this fact is  equivalent to the Brunn-Minkowski inequality for any convex bodies. 

If $K$ is origin symmetric, then  $-\mathcal{L}_K$ can be restricted to the space of even functions in $C^2(S^{n-1})$,
and $\lambda_{1,e}(-\mathcal{L}_K)>1$ holds for the smallest positive eigenvalue of this restricted operator because linear functions are odd.
Here the linear invariance yields that 
$\lambda_{1,e}\left(-\mathcal{L}_K\right)=\lambda_{1,e}\left(-\mathcal{L}_{\Phi K}\right)$ for $\Phi\in{\rm GL}(n)$.
A key result in Kolesnikov, Milman \cite{KoM22}  improves 
the estimate $\lambda_{1,e}(-\mathcal{L}_K)>1$ uniformly; more precisely,
$$
\lambda_{1,e}\left(-\mathcal{L}_K\right)\geq\frac{n-p_n}{n-1}
$$
for any $K\in \mathcal{K}_e^n$ with $C^2_+$ boundary where 
the explicit $p_n\in(0,1-\frac{c}{n(\log n)^{10}})$ is the same as in Theorem~\ref{pcloseto1}.
The connection to the $L_p$-Minkowski conjecture for fixed $p\in[0,1)$ is another key result in Kolesnikov, Milman \cite{KoM22}, as developed futher by Putterman \cite{Put21}; namely, 
\begin{equation}
\label{lambda1ep}
\lambda_{1,e}\left(-\mathcal{L}_K\right)\geq \frac{n-p}{n-1}
\end{equation}
is equivalent saying that \eqref{LpBMPuttermann} holds for any $C\in \mathcal{K}_e^n$. In particular, given
$p\in[0,1)$, the $L_p$-Minkowski conjecture follows if
\eqref{lambda1ep} holds for all $K\in \mathcal{K}_e^n$ with $C^\infty_+$ boundary, and the Logarithmic-Minkowski conjecture 
and \eqref{logBMconjeq} are equivalent saying that
\begin{equation}
\label{lambda1elog}
\lambda_{1,e}\left(-\mathcal{L}_K\right)\geq \frac{n}{n-1}
\end{equation}
 for all $K\in \mathcal{K}_e^n$ with $C^\infty_+$ boundary.
If $K_m\in \mathcal{K}_e^n$ with $C^2_+$ boundary tends to a cube, then 
$\lambda_{1,e}\left(-\mathcal{L}_{K_m}\right)$ tends to $\frac{n}{n-1}$ according to \cite{KoM22}; therefore, the Logarithmic-Minkowski conjecture states that cubes "minimize" $\lambda_{1,e}\left(-\mathcal{L}_K\right)$.
On the other hand, \cite{KoM22} calculates that $\lambda_{1,e}\left(-\mathcal{L}_{K}\right)=\frac{2n}{n-1}$
if $K$ is a centered Euclidean ball, and Milman \cite{Mila} verifies that centered ellipsoids maximize 
$\lambda_{1,e}\left(-\mathcal{L}_{K}\right)$ among all $K\in \mathcal{K}_e^n$ with $C^2_+$ boundary.

\section{Some variants of the $L_p$-Minkowski problem} 
\label{secVariants}

We note that Livshyts \cite{Liv19} considers a version of the Minkowski problem with a given measure on $\R^n$ acting as a weigh on the surface of the convex body. 

Considering the variation of the $i$th intrinsic volume of a convex body $K$ (or equivalently, variation of
$V(B^n,K;i)$) for $i=2,\ldots,n-1$ instead of the volume of $K$ leads to the so-called Christoffel-Minkowski problem, which asks to determine a convex body when its $(i-1)$th area
measure on $S^{n-1}$ is prescribed (see  Guan, Ma \cite{GuN17}, Guan, Xia \cite{GuX18}). 
We note that for $K\in\mathcal{K}^n$ with $C^2_+$ boundary, $S_K$ is then $(n-1)$th surface area measure, and the 
$j$th area measure is defined using the $j$th symmetric function of the principle radii of curvatures instead of the reciprocal of the Gaussian curvature. The $L_p$ Christoffel-Minkowski problem is discussed by 
Guan, Xia \cite{GuX18},
Hu, Ma, Shen \cite{HMS04} and Bryan, Ivaki, Scheuer \cite{BIJ19}  in the case $p>1$, and by
Bianchini, Colesanti, Pagnini, Roncoroni \cite{BCPR} in the case $p\in[0,1)$ where again, $p=1$ corresponds to the classical case.

The Minkowski problem on the sphere is solved by Guang, Li, Wang \cite{LGWb} (see \cite{LGWb} for related references, as well),
and in the hyperbolic space, partial results, also about the hyperbolic Christoffel-Minkowski problem, are obtained by
Gerhardt \cite{Ger}.

The Gaussian surface area measure of a $K\in\mathcal{K}^n$ is defined by
Huang, Xi and Zhao \cite{HXZ21}, and \cite{HXZ21} obtains significant results about the even Gaussian Minkowski problem.
These results are extended to the not necessarily even case by Feng, Liu, Xu \cite{FLX23} and 
Shibing Chen, Shengnan Hu, Weiru Liu, Yiming Zhao \cite{CHLZ}, and the 
$L_p$-Gaussian Minkowski problem is considered by Liu \cite{Liu22} and Feng, Hu, Xu \cite{FHX23}. Uniqueness of the solution is discussed in the works above and in Ivaki, Milman \cite{IvM}.

Next we discuss the $L_p$ dual Minkowski problem introduced by Lutwak, Yang, Zhang \cite{LYZ18}
that is a common generalization of the $L_p$-Minkowski problem and the Alexandrov problem. In order to define the dual curvature measures, let $K\in\mathcal{K}_{(o)}^n$.  
Recall that the radial function $\varrho_K(u)>0$ satisfies $\varrho_K(u)u\in\partial K$ for any $u\in S^{n-1}$. For a Borel set $\omega\subset S^{n-1}$, its $\HH$ measurable reverse radial Gauss image $\alpha^*(\omega)$
is the set of $u\in S^{n-1}$ such that some $v\in \omega$ is an exterior normal at $\varrho_K(u)u$ (see Huang, Lutwak, Yang, Zhang \cite{HLYZ16}). Now for any
$q\in \R$, \cite{HLYZ16} defines the $q$th dual curvature measure of the Borel set $\omega\subset S^{n-1}$ by
$$
\widetilde{C}_{K,q}(\omega)=\frac1n\int_{\alpha^*(\omega)}\varrho_K^n\,d\HH.
$$
In particular, $\widetilde{C}_{K,n}=V_K$ is the cone volume measure (discussed in Section~\ref{secLogMinkowski}), and  $n\widetilde{C}_{K,0}$ is the Alexandrov integral curvature measure of the polar $K^*$. The Monge-Amp\'ere equation corresponding to the $q$th dual Minkowski problem is
\begin{equation}
\label{MongeDualMink}
(\|\nabla h\|^2+h^2)^{\frac{q-n}2} \cdot h\det(\nabla^2h+h\,{\rm Id})= f.
\end{equation}
The Alexandrov  problem; namely, the charaterization of $\widetilde{C}_{K,0}$
has been solved by Aleksandrov \cite{Ale42,Ale96} (see also  Oliker \cite{Oli07} and Bertrand \cite{Ber16}
B\"or\"oczky, Lutwak, Yang, Zhang, Zhao \cite{BLYZ20}). For the $L_p$ version of the Alexandrov problem posed by Huang, Lutwak, Yang, Zhang \cite{HLYZ18}, see for example  
Zhao \cite{Zha19}, Li, Sheng, Ye, Yi \cite{LSYY22}
 Mui \cite{Mui22} and Wu, Wu,  Xiang \cite{WWX22}. If $q\neq 0,n$, then the following results are known:
\begin{itemize}
\item If $q<0$, then any Borel measure on $S^{n-1}$
not concentrated on a closed hemisphere is a $q$th dual Minkowski curvature measure according to 
Zhao \cite{Zha17} and Li, Sheng, Wang \cite{LSW20}.
\item If $0<q<n$, then an even Borel measure on $S^{n-1}$ is a $q$th dual Minkowski curvature measure if and only if 
$$
\mu(L\cap S^{n-1})<\mbox{$\frac{{\rm dim}\,L}q$}\cdot \mu(S^{n-1})
$$ 
for any proper linear subspace $L$ of of $\R^n$ according to B\"or\"oczky, Lutwak, Yang, Zhang, Zhao \cite{BLYZ19}
where one needs to add that $\mu$ is not concentrated onto a great subsphere if $q<1$.
\item If $q\geq n+1$ and $K\in\mathcal{K}^n$ is origin symmetric, then Henk, Pollehn \cite{HeP18} prove
$$
\widetilde{C}_{K,q}(L\cap S^{n-1})<\mbox{$\frac{q-n+{\rm dim}\,L}q$}\cdot \widetilde{C}_{K,q}(S^{n-1})
$$ 
\item If $q>0$ and $n=2$, then  \eqref{MongeDualMink} has a solution for any measurable $f$ provided $\frac1{c}<f<c$ for a $c>1$ according to Chen, Li \cite{ChL18}.
\end{itemize}
In particular, it is an intriguing open problem to characterize an even $q$th dual Minkowski curvature measure on $S^{n-1}$ if $q>n$.

For $p,q\in\R$, Lutwak, Yang, Zhang \cite{LYZ18} defines the $q$th $L_p$ dual Minkowski curvature measure on $S^{n-1}$ by $d\widetilde{C}_{K,p,q}=h_K^{-p}d\widetilde{C}_{K,q}$,
and hence $\widetilde{C}_{K,0,q}=\widetilde{C}_{K,q}$ and $\widetilde{C}_{K,p,n}=\frac1n\,S_{K,p}$.
Given a Borel measure $\mu$ on $S^{n-1}$, the simplest version of the $q$th the $L_p$ dual Minkowski problem
asks for a $K\in\mathcal{K}_{(o)}^n$ with $\mu=\widetilde{C}_{K,p,q}$, and the correspong Monge-Amp\'ere equation is
\begin{equation}
\label{LpdualMinkowski}
 h^{1-p}\det(\nabla^2h+h\,{\rm Id})= (\|\nabla h\|^2+h^2)^{\frac{n-q}2} \cdot f
\end{equation}
Improving on \cite{BoF19} and Huang, Zhao \cite{HuZ18}, Chen, Li \cite{ChL21} and
Lu, Pu \cite{LuP21}
 if $p>0$ and $q\neq p$, then any Borel measure
not concentrated on a closed hemisphere is a $q$th $L_p$ dual Minkowski curvature measure (more precisely, if $p\leq q$, then some modification of the Monge-Amp\'ere equation might be needed). 
Huang, Zhao \cite{HuZ18} proved the same for $p,q<0$ and $p\neq q$ within the category of even measures.
See also Guang, Li, Wang \cite{LGW23a} for a flow approach when $p<0$ and $q>n$ under regularity assumptions.

Uniqueness of the solution of the $q$th $L_p$ dual Minkowski problem \eqref{LpdualMinkowski} is thouroughly investigated by Li, Liu, Lu \cite{LLL22}. The case when $n=2$ and $f$ is a constant function has been completely clarified by
Li,  Wan \cite{LiW}.

Another important related variant of the dual Minkowski problem is the so-called Chord Minkowski Problem
({\it cf.} Lutwak, Xi, Yang, Zhang \cite{LXYZ}) and its $L_p$ version by  Xi, Yang, Zhang, Zhao \cite{XYZZ23} for $p>0$, and 
by Li \cite{YLia,YLib} for $p<0$, see also  Guo,  Xi, Zhao \cite{GXZ} and 
Xi, Yang, Zhang, Zhao \cite{XYZZ23}.
In addition, the
"Affine dual Minkowski problem" is proposed by Cai, Leng, Wu,  Xi \cite{CLWX}.

Cordero-Erausquin, Klartag \cite{CEK15} generalized and solved the classical Minkowski problem to the space of log-concave functions on $\R^n$. Recently, Fang, Xing, Ye  \cite{FXY22} considered the $L_p$-Minkowski problem for log-concave functions,
and
Huang, Liu, Xi, Zhao \cite{HLXZ} and Fang, Ye, Zhang, Zhao \cite{FYZZ} managed to extend the dual Minkowski problem to
 log-concave functions on $\R^n$ for $q>0$.

Starting with Haberl, Lutwak, Yang, Zhang \cite{HLYZ10}, Orlicz versions of the $L_p$-Minkowski problem have been intensively investigated; namely, the function $t\mapsto t^{1-p}$ in \eqref{MongeSKp} is replaced by certain
$\varphi:(0,\infty)\to(0,\infty)$, and hence \eqref{MongeSKp} is replaced by
$$
\varphi(h)\det(\nabla^2 h+h\,{\rm Id})= f
$$
where $f$ is a given non-negative function on $S^{n-1}$.
Typically, the solution is only up to a constant factor; namely, there exists some $c>0$ such that 
$\varphi(h)\det(\nabla^2 h+h\,{\rm Id})= c\cdot f$.
The known existence results about the $L_p$-Minkowski problem for $p>-n$ have been generalized to the
 Orlicz $L_p$-Minkowski problem where $\varphi(t)$ replaces $t^{1-p}$  by Huang, He \cite{HuH12}
if $p>1$ (see also Xie \cite{Xie22}), by Jian, Lu \cite{JiL19}  if $p\in(0,1)$, and by Bianchi, B\"or\"oczky, Colesanti \cite{BBC19}
if $p\in(-n,0)$. 

Orlicz versions of these Monge-Amp\`ere equations have been 
considered by  Li, Sheng,  Ye, Yi \cite{LSYY22}
and Feng, Hu, Liu \cite{FHL22} in the case of the Alexandrov problem, by Xing, Ye \cite{XiZ20} in the case
of the dual Minkowski problem, and
Gardner, Hug, Weil, Xing, Ye \cite{GHWXY19,GHXY20},
Xing, Ye, Zhu \cite{XYZ22} and 
Liu, Lu \cite{LiL20} in the case
of the $L_p$ dual Minkowski problem in general.\\

\noindent{\bf Acknowledgements: } I would like to thank for illuminating discussions with Emanuel Milman, Erwin Lutwak, Gaoyong Zhang, Shibing Chen, Qi-Rui Li, Christos Saroglou and Martin Henk during the preparation of the manuscript. Also grateful for the hospitality of ETH Z\"urich during the completion of the survey.

\end{document}